\documentclass[reqno,11pt]{amsart}

\usepackage[cmtip,arrow]{xy}  % commuting diagrams
\usepackage{pb-diagram,pb-xy}  % commuting diagrams
\usepackage{graphicx}    % standard LaTeX graphics tool when including figure files
\usepackage{amsfonts}
\usepackage{amssymb}
\usepackage[latin1]{inputenc}

\newtheorem{thm}{Theorem}[section]
\newtheorem{cor}[thm]{Corollary}
\newtheorem{lemma}[thm]{Lemma}

\newtheorem{prop}[thm]{Proposition}
\theoremstyle{definition}
\newtheorem{defn}[thm]{Definition}
\theoremstyle{remark}

\newtheorem{remark}[thm]{Remark}
\newtheorem{example}[thm]{Example}

\newcommand{\PB}{\left\{\cdot\,,\cdot\right\}}

\newcommand{\Pb}[1]{\left\{\cdot\,,#1\right\}}
\newcommand{\pb}[1]{\left\{#1\right\}}

\newcommand{\lb}[1]{\[#1\]}

\renewcommand{\(}{\left(}
\renewcommand{\)}{\right)}
\renewcommand{\[}{\left[}
\renewcommand{\]}{\right]}

\newcommand{\set}[1]{\left\{#1\right\}}

\newcommand{\leqs}{\leqslant}
\newcommand{\geqs}{\geqslant}

\newcommand{\pp}[2]{\frac{\partial#1}{\partial#2}}
\newcommand{\p}{\partial}
\newcommand{\we}{\wedge}
\newcommand{\Rk}{\hbox{Rk\,}}
\renewcommand{\l}{\lambda}

\newcommand{\F}{{\bf F}}
\newcommand{\G}{{\bf G}}

\newcommand{\X}{\mathcal X}
\newcommand{\Y}{\mathcal Y}

\newcommand{\D}{\mathcal D}
\renewcommand{\L}{\mathcal L}

\newcommand{\U}{\mathcal U}
\newcommand{\Vf}{\mathcal V}

\newcommand{\step}[1]{\underline{Step #1}. }
\newcommand{\rp}{{r'}}
\newcommand{\V}{{\mathcal V}}
\newcommand{\R}{\mathbf{R}}
\newcommand{\Z}{\mathbf{Z}}

\newcommand{\T}{\mathbf{T}}
\newcommand{\diff}{{\rm d }}
\newcommand{\FF}{\mathcal F}
\newcommand{\GG}{\mathcal G}

\newif\ifprivate
\privatefalse

 \numberwithin{equation}{section}

\def\???{\ifprivate {\bf {???}} \marginpar{{\Huge {\bf ?}}}\else \fi}
\numberwithin{equation}{section}

\begin{document}
\nocite{*}

\title[Action-angle coordinates]{Action-angle coordinates for integrable systems on Poisson manifolds}

\author{Camille Laurent-Gengoux}\address{Camille Laurent-Gengoux, Laboratoire de Math\'ematiques et Applications,
         UMR 6086 CNRS, Universit\'e de Poitiers, Boulevard Marie et Pierre CURIE, BP 30179, 86962 Futuroscope
         Chasseneuil Cedex, France}\email{camille.laurent@math.univ-poitiers.fr}
\author{Eva Miranda${}^1$} \address{Eva Miranda, Departament de Matemàtiques, Universitat Autònoma de Bar\-celona, E-08193
         Bellaterra, Spain}\email{miranda@mat.uab.cat}
\thanks{${}^1$Research supported by a Juan de la Cierva contract and partially supported by the DGICYT/FEDER, project number
         MTM2006-04353 (Geometr\'\i a Hiperb\'olica y Geometr\'\i a Simpl\'ectica).}
\author{Pol Vanhaecke${}^2$} \address{Pol Vanhaecke, Laboratoire de Math\'ematiques et Applications,
         UMR 6086 du CNRS, Universit\'e de Poitiers, Boulevard Marie et Pierre CURIE, BP 30179, 86962 Futuroscope
         Chasseneuil Cedex, France}\email{pol.vanhaecke@math.univ-poitiers.fr}
\thanks{${}^2$Partially supported by a European Science Foundation grant (MISGAM), a Marie Curie grant (ENIGMA) and
         an ANR grant (GIMP)}

\date{\today}


% ----------------------------------------------------------------

\begin{abstract}

We prove the action-angle theorem in the general, and most natural, context of integrable systems on Poisson
manifolds, thereby generalizing the classical proof, which is given in the context of symplectic manifolds. The
topological part of the proof parallels the proof of the symplectic case, but the rest of the proof is quite
different, since we are naturally led to using the calculus of polyvector fields, rather than differential forms;
in particular, we use in the end a Poisson version of the classical Carathéodory-Jacobi-Lie theorem, which we also
prove. At the end of the article, we generalize the action-angle theorem to the setting of non-commutative
integrable systems on Poisson manifolds.

\end{abstract}

\maketitle

\tableofcontents

% -*- mode: latex; tex-main-file: "art.tex"; -*-
%
\section{Introduction}
The action-angle theorem is one of the basic theorems in the theory of integrable systems. In this paper we prove
this theorem in the general, and most natural, context of integrable systems on Poisson manifolds.

\smallskip

We recall that a Poisson manifold $(M,\Pi)$ is a smooth manifold $M$ on which there is given a bivector field
$\Pi$, with the property that the bracket on $C^\infty(M)$, defined for arbitrary smooth functions $f$ and $g$ on
$M$ by
\begin{equation*}
  \pb{f,g}:=\Pi(\diff f,\diff g)
\end{equation*}%
is a Lie bracket, i.e., it satisfies the Jacobi identity. On a Poisson manifold $(M,\Pi)$, the Hamiltonian
operator, which assigns to a function on $M$ a vector field on $M$, is defined naturally by contracting the
bivector field with the function (the ``Hamiltonian''): for $h\in C^\infty(M)$ its Hamiltonian vector field is
defined by
\begin{equation}\label{for:intro_ham}
  \X_h:=\Pb h=-\imath_{\diff h}\Pi.
\end{equation}%
Two important consequences of the Jacobi identity for $\PB$ are that the (generalized distribution) on $M$, defined
by the Hamiltonian vector fields $\X_h$ is integrable, and that the Hamiltonian vector fields which are associated
to Poisson commuting functions (usually called functions in involution) are commuting vector fields.  The main
examples of Poisson manifolds are symplectic manifolds and the dual of a (finite-dimensional) Lie algebra, but
there are many other examples, which come up naturally in deformation theory, the theory of $R$-brackets,
Lie-Poisson groups, and so on. Poisson's original bracket on $C^{\infty}(\R^{2r})$, given for smooth functions $f$
and $g$ by
\begin{equation}\label{poisson_form_intro}
  \pb{f,g}:=\sum_{i=1}^r\left(\frac{\p f}{\p q_i}\frac{\p g}{\p p_i}-\frac{\p g}{\p q_i}\frac{\p f}{\p p_i}\right),
\end{equation}
is still today of fundamental importance in classical and quantum mechanics, and in other areas of mathematical
physics. Many examples of integrable Hamiltonian systems are known in the context of Poisson manifolds which are
not symplectic. For instance the Kepler problem \cite{Whi}, Toda lattices \cite{adlermoerbekevanhaecke2004} and
the Gelfand-Cetlin systems \cite{guilleminandsternberg,Giacobbe}.

\smallskip

One of the main uses of the Poisson bracket is the integration of
Hamilton's equations, which are the equations of motion which
describe a classical mechanical system on the phase space
$\R^{2r}\simeq T^*\R^r$, defined by a Hamiltonian $h$ (the energy,
viewed as a function on phase space); their solutions are the
integral curves of the Hamiltonian vector field $\X_h$, defined
by~(\ref{for:intro_ham}) with respect to the Poisson bracket
(\ref{poisson_form_intro}). The fundamental Liouville theorem states
that it suffices to have $r$ independent functions in involution
$(f_1=h,f_2,\dots,f_r)$ to quite explicitly (i.e., by quadratures)
integrate the equations of motion for generic initial conditions.
Moreover, assuming that the so-called invariant manifolds, which are
the (generic) submanifolds traced out by the $n$ commuting vector
fields $\X_{f_i}$, are compact, they are (diffeomorphic to) tori
$\T^r=\R^r/\Lambda$, where $\Lambda$ is a lattice in $\R^r$; on
these tori, which are known as Liouville tori, the flow of each of
the vector fields $\X_{f_i}$ is linear, so that the solutions of
Hamilton's equations are quasi-periodic. The classical action-angle
theorem goes one step further: under the above topological
assumption, there exist on a neighbourhood $U$ of every Liouville
torus functions $\sigma_1,\dots,\sigma_r$ and $\R/\Z$-valued
functions $\theta_1,\dots,\theta_r$, having the following
properties:
\begin{enumerate}
  \item[(1)] The map $\Phi$, defined by $\Phi:=(\theta_1,\dots,\theta_r,\sigma_1,\dots,\sigma_r)$ is a
    diffeomorphism from $U$ onto the product $\T^r\times B^r$, where $B^r$ is an $r$-dimensional ball;
  \item[(2)] $\Phi$ is a canonical map: in terms of $\theta_1,\dots,\theta_r,\sigma_1,\dots,\sigma_r$ the Poisson
    structure takes the same form as in (\ref{poisson_form_intro}) (upon replacing $q_i$ by $\theta_i$ and $p_i$ by
    $\sigma_i$);
  \item[(3)] Under $\Phi$, the Liouville tori in $U$ correspond to the fibers of the natural projection $\T^r\times
    B^r\to B^r$.
\end{enumerate}
The proof of this theorem goes back to Mineur \cite{Mineur-AA1935,Mineur-AA1936,Mineur-AA1937}. A proof in the case
of a Liouville integrable system on a symplectic manifold was given by Arnold \cite{arnold}; see also
\cite{batesandcushman,Duistermaat,libermannmarle}.
As established in \cite{guilleminandsternberg}, action-angle coordinates also appear naturally in geometric
quantization, for, when an integrable system is interpreted as a polarization, action-angle coordinates determine
the so-called Bohr-Sommerfeld leaves: the latter are in particular explicitely described for the Gelfand-Cetlin
system in \cite{guilleminandsternberg}.

\smallskip

In the context of Poisson manifolds, the Liouville theorem still holds, up to two adaptations: one needs to take into
account the Casimirs (functions whose Hamiltonian vector field are zero) and the singularities of the Poisson
structure (the points where the rank of the bivector field drops); for a precise statement and a proof, see
\cite[Ch.\ 4.3]{adlermoerbekevanhaecke2004}. As we show in this paper, the action-angle theorem takes in the case
of Poisson manifolds the following form\footnote{An equivalent statement, without proof, was given in
\cite[Ch.\ 4.3]{adlermoerbekevanhaecke2004}.}

\begin{thm}\label{thm:action-angle_intro}
  Let $(M,\Pi)$ be a Poisson manifold of dimension $n$ and (maximal) rank~$2r$. Suppose that $\F=(f_1,\dots,f_s)$
  is an integrable system on $(M,\Pi)$, i.e., $r+s=n$ and the components of $\F$ are independent and in
  involution. Suppose that $m\in M$ is a point such that
\begin{enumerate}
  \item[(1)] $\diff_mf_1\we\dots\we\diff_mf_s\neq0$;
  \item[(2)] The rank of $\Pi$ at $m$ is $2r$;
  \item[(3)] The integral manifold $\FF_m$ of $\X_{f_1},\dots,\X_{f_s}$, passing through $m$, is compact.
\end{enumerate}
  Then there exists ${\R} $-valued smooth functions $(\sigma_1,\dots, \sigma_{s})$ and $ {\R}/{\Z}$-valued smooth
  functions $({\theta_1},\dots,{\theta_r})$, defined in a neighborhood $U$ of $\FF_m$ such that
  \begin{enumerate}
    \item The functions $(\theta_1,\dots,\theta_r,\sigma_1,\dots,\sigma_{s})$ define an isomorphism
      $U\simeq\T^r\times B^{s}$;
    \item The Poisson structure can be written in terms of these coordinates as
      $$
        \Pi=\sum_{i=1}^r\pp{}{\theta_i}\we\pp{}{\sigma_i},
      $$
      in particular the functions $\sigma_{r+1},\dots,\sigma_{s}$ are Casimirs of $\Pi$ (restricted to $U$);
    \item The leaves of the surjective submersion $\F=(f_1,\dots,f_{s})$ are given by the projection onto the
          second component $\T^r\times B^{s}$, in particular, the functions $\sigma_1,\dots,\sigma_{s}$ depend on the
          functions $f_1,\dots,f_{s}$ only.
  \end{enumerate}
\end{thm}
The functions $\theta_1,\dots,\theta_r$ are called {angle coordinates}, the functions $\sigma_1,$
$\dots,\sigma_r$ are called {action coordinates} and the remaining functions
$\sigma_{r+1},\dots,\sigma_{s}$ are called {transverse coordinates}.

\smallskip

Our proof of theorem \ref{thm:action-angle_intro}, consists of several conceptually different steps, which are in
$1$-$1$ correspondence with the (a) topological, (b) group theoretical, (c) geometrical and (d) analytical aspects
of the construction of the coordinates. It parallels Duistermaat's proof, which deals with the symplectic case
\cite{Duistermaat}; while (a) and (b) are direct generalizations of his proof, (c) and (d) are however not.

\smallskip

(a) The topological part of the proof amounts to showing that in the neighborhood of the invariant manifold
$\FF_m,$ we have locally trivial torus fibration (Paragraph \ref{par:fol_by_tori}). Once we have shown that the
compact invariant manifolds are the connected components of the fibers of a submersive map (the map $\F$,
restricted to some open subset), the proof of this part is similar as in the symplectic case.

\smallskip

(b) The (commuting) Hamiltonian vector fields are tangent to the
tori of this fibration; integrating them we get an induced torus
action (action by~$\T^r$) on each of these tori, but in general
these actions cannot be combined into a single torus action. Taking
appropriate linear combinations of the vector fields, using
$\F$-basic functions as coefficients, by a procedure called
``uniformization of the periods'', one constructs new vector fields
$Y_1,\dots,Y_r$ which are tangent to the fibration, and which now
integrate into a single torus action. This is the content of step 1
in the proof of proposition \ref{prp:actions1}. This step, which is
an application of the implicit function theorem, is identical as in
the symplectic case.

\smallskip

(c) The newly constructed vector fields $Y_i$ are the fundamental vector fields of a torus action. We first show
that they are Poisson vector fields, i.e., that they preserve the Poisson structure (step 2 in the proof of
proposition~\ref{prp:actions1}). The key (and non-trivial) point of the proof is the periodicity of the vector
fields $Y_i$. We then prove (in step 3 of the proposition) the stronger statement that these vector fields are
Hamiltonian vector fields, at least locally, by constructing quite explicitly their Hamiltonians, which will in the
end play the role of action coordinates.

\smallskip

(d) In the last step (theorem \ref{thm:action-angle}), we use the Carathéodory-Jacobi-Lie theorem for Poisson
manifolds to construct on the one hand coordinates which are conjugate to the action coordinates (angle coordinates
and transverse coordinates) and on the other hand to extend these coordinates to a neighborhood of the Liouville
torus $\FF_m$. The Carathéodory-Jacobi-Lie theorem for Poisson manifolds, to which Section
\ref{sec:darboux-caratheodory} is entirely devoted, provides a set of canonical local coordinates for a Poisson
structure $\Pi$, containing a given set $p_1,\dots,p_r$ of functions in involution. It generalizes both the
classical Carathéodory-Jacobi-Lie theorem for symplectic manifolds \cite[Th.\ 13.4.1]{libermannmarle} and
Weinstein's splitting theorem \cite[Th.\ 2.1]{weinstein}. We are convinced that this theorem, which is new, has
other interesting applications, as in the study of local forms and stability of integrable systems.

\medskip

The action-angle theorem has been proven by \cite{Nekhroshev} in the general context of non-commutative integrable
systems on a symplectic manifold (see the Appendix for a comparison between this notion and some closely related
notions of integrability). Roughly speaking, a non-commutative integrable system has more constants of motion than
a Liouville integrable system, accounting for linear motion on smaller tori, but not all these functions are in
involution. This notion has a natural definition in the case of Poisson manifolds, proposed here (definition
\ref{def:non-com}); it generalizes both the notion of Liouville integrability on a Poisson manifold and the notion
of non-commutative integrability on a symplectic manifold. We show in Section \ref{sec:non-com} that our proof can
be adapted (i.e., generalized) to provide a proof of the action-angle theorem in this very general context.

\medskip

The structure of the paper is as follows. We state and prove the Carathéo\-dory-Jacobi-Lie theorem for Poisson
manifolds in Paragraph \ref{par:darboux-caratheodory} and we give in Paragraph \ref{par:example} a counterexample
which shows that a mild generalization of the latter theorem does not hold in general. The action-angle theorem for
Liouville integrable systems on Poisson manifolds is given in Section \ref{sec:aa_coordinates}. We show in Section
\ref{sec:non-com} how this theorem can be adapted to the more general case of non-commutative integrable systems on
Poisson manifolds. The appendix to the paper is devoted to the geometrical formulation of the notion of a
non-commutative integrable system on a Poisson manifold.

\smallskip

In this paper, all manifolds and objects considered on them are smooth and we write $\pb{f,g}$ for $\Pi(\diff
f,\diff g)$.

% -*- mode: latex; tex-main-file: "art.tex"; -*-
%
\section{The Carathéodory-Jacobi-Lie theorem for Poisson manifolds}
\label{sec:darboux-caratheodory}

In this section we prove a natural generalization of the classical Carathéo\-dory-Jacobi-Lie theorem \cite[Th.\
 13.4.1]{libermannmarle} for an arbitrary Poisson manifold $(M,\Pi)$. It provides a set of canonical local
coordinates for the Poisson structure $\Pi$, which contains a given set $p_1,\dots,p_r$ of functions in involution
(i.e., functions which pairwise commute for the Poisson bracket), whose Hamiltonian vector fields are assumed to be
independent at a point $m\in M$ (theorem \ref{thm:localsplitting}). This result, which is interesting in its own
right, will be used in our proof of the action-angle theorem. We show in Paragraph \ref{par:example} by giving a
counterexample that canonical coordinates containing a given set of functions in involution may fail to exist as
soon as the Hamiltonian vector fields $\X_{p_1},\dots,\X_{p_r}$ are dependent at $m$, even if they are independent
at all other points in a neighborhood of $m$.

\subsection{The theorem}
\label{par:darboux-caratheodory}

The main result of this section is the following theorem.

\begin{thm} \label{thm:localsplitting}
  Let $m$ be a point of a Poisson manifold $(M,\Pi)$ of dimension~$n$. Let $p_1,\dots,p_r$ be $r$ functions in
  involution, defined on a neighborhood of $m$, which vanish at $m$ and whose Hamiltonian vector fields are
  linearly independent at $m$.  There exist, on a neighborhood $U$ of $m$, functions
  $q_1,\dots,q_r,z_1,\dots,z_{n-2r}$, such that
  \begin{enumerate}
    \item The $n$ functions $(p_1,q_1,\dots,p_r,q_r, z_1,\dots,z_{n-2r})$ form a system of coordinates on $U$,
    centered at $m$;
    \item The Poisson structure $\Pi$ is given on $U$ by
    \begin{equation}\label{eq:thm_split}
      \Pi=\sum_{i=1}^r\pp{}{q_i}\we\pp{}{p_i}+\sum_{i,j=1}^{n-2r} g_{ij}(z)\pp{}{z_i}\we\pp{}{z_j},
    \end{equation}
    where each function $g_{ij}(z)$ is a smooth function on $U$ and is independent of $p_1,\dots,p_r,
    q_1,\dots,q_r$.
  \end{enumerate}
  The rank of $\Pi$ at $m$ is $2r$ if and only if all the functions $g_{ij}(z)$ vanish for $z=0$.
%  Such local coordinates $(p_1,q_1,\dots,p_r,q_r, z_1,\dots,z_{n-2r}) $ are called \emph{coordinates conjugated to
%  $p_1,\dots,p_r$}.
%
\end{thm}
\begin{proof}
We show the first part of the theorem by induction on $r$. For $r=0$, every system of coordinates $z_1,\dots,z_n$,
centered at $m$, does the job.  Assume that the result holds true for every point in every Poisson manifold and
every $(r-1)$-tuple of functions as above, with $r\geqs1$. We prove it for~$r$. To do this, we consider an
arbitrary point $m$ in an $n$-dimensional Poisson manifold $(M,\Pi)$, and we assume that we are given functions in
involution $p_1,\dots,p_r$, defined on a neighborhood of $m$, which vanish at $m$, and whose Hamiltonian vector
fields are linearly independent at $m$. On a neighbourhood of $m$, the distribution
$\D:=\langle{\X}_{p_1},\dots,{\X}_{p_{r}} \rangle$ has constant rank $r$ and is an involutive distribution because
$[\X_{p_i},\X_{p_j}]= -\X_{\pb{p_i,p_j}}=0$. By the Frobenius theorem, there exist local coordinates $g_1,\dots,
g_{n}$, centered at $m$, such that $\X_{p_i}=\pp{}{g_i}$ for $i=1,\dots,r$, on a neighbourhood of $m$. Setting
$q_{r}:=g_{r}$ we have
  \begin{equation} \label{eq:q1}
    {\X}_{q_r}[p_i]=-{\X}_{p_i}[q_r]=-\delta_{i,r}, \qquad i=1,\dots,r,
  \end{equation}
in particular (1) the $r+1$ vectors $\diff_m p_1,\dots,\diff_m p_{r}$ and $\diff_m q_{r}$ of $T^*_mM$ are linearly
independent, and (2) the vector fields $\X_{q_r}$ and $\X_{p_r}$ are independent at~$m$. It follows that a
distribution $\D'$ (of rank $2$) is defined by $\X_{q_r}$ and~$\X_{p_r}$. It is an integrable distribution because
$[{\X}_{q_{r}},{\X}_{p_{r}}]=-{\X}_{\pb{q_{r},p_{r}}}=0.$ Applied to $\D'$, the Frobenius theorem yields the
existence of local coordinates $v_1,\dots,v_{n}$, centered at $m$, such that
\begin{equation}\label{eq:r_der}
  {\X}_{p_r}=\pp{}{v_{n-1}}\qquad\hbox{and}\qquad{\X}_{q_r}=\pp{}{v_n}.
\end{equation}%
Since the differentials $ \diff_m v_1,\dots,\diff_m v_{n-2}$ vanish on ${\X}_{p_r}(m)$ and on ${\X}_{q_r}(m)$, it
follows that $(\diff_mv_1,\dots,\diff_m v_{n-2},\diff_m p_r,\diff_m q_r) $ is a basis of $T^*_mM$.  Therefore, the
$n$ functions $(v_1,\dots,v_{n-2},p_{r},q_{r})$ form a system of local coordinates, centered at $m$.  It follows
from (\ref{eq:r_der}) that the Poisson structure takes in terms of these coordinates the following form:
\begin{equation*}
  \Pi=\pp{}{q_r}\we\pp{}{p_r}+\sum_{i,j=1}^{n-2}h_{ij}(v_1,\dots,v_{n-2},p_r,q_r)\pp{}{v_i}\we\pp{}{v_j}.
\end{equation*}%
The Jacobi identity, applied to the triplets $(p_{r},v_i,v_j)$ and $(q_{r},v_i,v_j) $, implies that the functions
$h_{ij}$ do not depend on the variables $p_{r},q_{r}$, so that
\begin{equation} \label{eq:local1}
  \Pi=\pp{}{q_r}\we\pp{}{p_r}+\sum_{i,j=1}^{n-2}h_{ij}(v_1,\dots,v_{n-2})\pp{}{v_i}\we\pp{}{v_j},
\end{equation}
which means that $\Pi$ is, in a neighborhood of $m$, the product of a symplectic structure (on a neighborhood of
the origin in $\R^2$) and a Poisson structure (on a neighborhood of the origin in $\R^{n-2}$). In order to apply
the recursion hypothesis, we need to show in case $r-1>0$ that $p_1,\dots,p_{r-1}$ depend only on the coordinates
$v_1,\dots,v_{n-2}$, i.e., are independent of $p_r$ and $q_r$,
\begin{equation}\label{eq:indep_r}
   \pp{p_i}{p_r} = 0=\pp{p_i}{q_r}\qquad i=1,\dots,r-1.
\end{equation}%
Both equalities in (\ref{eq:indep_r}) follow from the fact that $p_i$ is in involution with $p_r$ and $q_r$, for
$i=1,\dots, r-1$, combined with (\ref{eq:local1}):
\begin{equation*}
  0=\pb{p_i,p_r}=\pp{p_i}{q_r},\qquad 0=\pb{p_i,q_r}=-\pp{p_i}{p_r}.
\end{equation*}%
We may now apply the recursion hypothesis on the second term in (\ref{eq:local1}), together with the functions
$p_1,\dots, p_{r-1}$. It leads to a system of local coordinates $(p_1,q_1,\dots,p_r,q_r, z_1,\dots,z_{n-2r})$ in
which $\Pi$ is given by (\ref{eq:thm_split}). This shows the first part of the theorem. The second part of the
theorem is an easy consequence of (\ref{eq:thm_split}), since it implies that the rank of $\Pi$ at $m$ is $2r$ plus
the rank of the second term in the right hand side of (\ref{eq:thm_split}), at $z=0$.
\end{proof}

\begin{remark}
  The classical Carathéodory-Jacobi-Lie theorem corresponds to the case $\dim M=2r$. Then $\Pi$ is the Poisson
  structure associated to a symplectic structure, in the neighborhood of $m$. Theorem \ref{thm:localsplitting} then
  says that $\Pi$ can be written in the simple form
  \begin{equation}\label{eq:DC}
    \Pi=\sum_{i=1}^r\pp{}{q_i}\we\pp{}{p_i},
  \end{equation}
  where we recall that the (involutive) set of functions $p_1,\dots,p_r$ is prescribed.
\end{remark}

\begin{remark}
  Theorem \ref{thm:localsplitting} and their proof, as they are stated, do not yield the existence of the
  involutive set of functions $p_1,\dots,p_r$, a fact which is plain in Weinstein's splitting theorem. However, if
  we forget in our proof that these functions are prescribed, we can easily adapt the induction hypotheses, adding
  the existence of $r$ such functions, when the rank of the Poisson structure at $m$ is at least $2r$. In this
  sense, our theorem is an amplification of Weinstein's splitting theorem.
\end{remark}

\begin{remark}
  Theorem \ref{thm:localsplitting} holds true for holomorphic Poisson manifolds; the local coordinates are in this
  case holomorphic coordinates and the functions $g_{ij}(z)$ are holomorphic functions, independent from
  $p_1,\dots,p_r, q_1,\dots,q_r$. Up to these substitutions, the given proof is valid word by word.
\end{remark}

\subsection{A counterexample}\label{par:example}
If we denote in theorem \ref{thm:localsplitting} the rank of $\Pi$ at $m$ by $2\rp$, then $2\rp\geqs 2r$, because
the involutive set of functions $p_1,\dots,p_r$ define a totally isotropic foliation in a neighborhood of $m$. It
means that, if $2\rp<2r$ and one is given independent functions in involution $p_1,\dots,p_r$, then their
Hamiltonian vector fields ${\X}_{p_1},\dots,{\X}_{p_r}$ are \emph{dependent} at $m$. In the extremal case in which
$\dim\left\langle{\X}_{p_1}(m),\dots,{\X}_{p_r}(m)\right\rangle=\rp$ one has\footnote{Possibly up to a relabelling
of the $p_i$, so that $\dim\left\langle{\X}_{p_1}(m),\dots,{\X}_{p_{\rp}}(m)\right\rangle=\rp$.}, according to
theorem \ref{thm:localsplitting}, that there exist functions $q_1,\dots,q_{\rp}$ and $z_1,\dots,z_{n-2\rp}$ such
that $\Pi$ takes the form
$$
  \Pi=\sum_{i=1}^{\rp}\pp{}{q_i}\we\pp{}{p_i}+\sum_{k,l=1}^{n-2\rp}\phi_{k,l}(z_1,\dots,z_{n-2r'})\pp{}{z_k}\we\pp{}{z_l}.
$$
A natural question is whether $r-\rp$ of the functions $z_i$ can be chosen as $p_{\rp+1},\dots,p_r$, or, more
generally, as functions which depend only on $p_1,\dots,p_r$. We show in the following (counter) example that this
is not possible, in general.

\begin{example}
On $\R^4$, with coordinates $f_1,f_2,g_1,g_2$, consider the bivector field, given by
\begin{equation}\label{eq:Pi_ori}
  \Pi=\pp{}{g_1}\we\pp{}{f_1}+\chi(g_2)\pp{}{g_2}\we\pp{}{f_2}+\psi(g_2)\pp{}{g_1}\we\pp{}{f_2},
\end{equation}
where $\chi(g_2)$ and $\psi(g_2) $ are smooth functions that depend only on $g_2$, and which vanish for $g_2=0$, so
that the rank of $\Pi$ at the origin is $2$.  A direct computation shows that this bivector field is a Poisson
bivector field and that $f_1$ and~$f_2$ are in involution. We show that for some choice of $\chi$ and $\psi$ there
exists no system of coordinates $p_1,q_1,z_1,z_2$, centered at $0$, with $p_1,z_1$ depending only on $f_1$ and
$f_2$, such that
\begin{equation}\label{eq:poisson_fic}
  \Pi =\pp{}{q_1}\we\pp{}{p_1}+\phi(z_1,z_2)\pp{}{z_1}\we\pp{}{z_2}.
\end{equation}
To do this, let us assume that such a system of coordinates exists. Taking the Poisson bracket of $p_1=p_1(f_1,f_2)$
and $z_1=z_1(f_1,f_2)$ with $q_1$ yields, in view of (\ref{eq:poisson_fic}),
\begin{eqnarray}\label{eq:lins}
  1&=&\pb{q_1,p_1}=\displaystyle\pp{p_1}{f_1}\pb{q_1,f_1}+\pp{p_1}{f_2}\pb{q_1,f_2},\nonumber\\
  0&=&\pb{q_1,z_1}=\displaystyle\pp{z_1}{f_1}\pb{q_1,f_1}+\pp{z_1}{f_2}\pb{q_1,f_2}.
\end{eqnarray}
Let $N$ denote the locus defined by $f_1=f_2=0$, which is a smooth surface in a neighborhood of the origin. Let $q$
denote the restriction of $q_1$ to $N$. Since $\X_{f_1}$ and $\X_{f_2}$ are tangent to $N$,
$\X_{f_i}[q]=\pb{q_1,f_i}_{\vert_N}$, so that (\ref{eq:lins}), restricted to $N$, becomes
\begin{eqnarray}\label{eq:ders}
  1&=&\displaystyle \l_1{\X}_{f_1}[q]+\l_2{\X}_{f_2}[q],\nonumber\\
  0&=&\displaystyle \l_3{\X}_{f_1}[q]+\l_4{\X}_{f_2}[q],
\end{eqnarray}
where $\l_1,\dots,\l_4$ are constants (because $p_1,z_1$ depend only on $f_1,f_2$), and satisfy
$\l_1\l_4-\l_2\l_3\neq0$, since $p_1$ and $z_1$ are part of a coordinate system centered at the origin. It follows
that
\begin{equation}\label{eq:cond_q}
  {\X}_{f_1}[q]=c_1  \mbox{ and }  {\X}_{f_2}[q] = c_2,
\end{equation}
where $c_1$ and $c_2$ are constants, which cannot be both equal to zero, in view of (\ref{eq:ders}). Writing
$\X_{f_1}$ and $\X_{f_2}$ in terms of the original variables, using (\ref{eq:Pi_ori}), we find that $q=q(g_1,g_2)$
must satisfy
$$
  \pp q{g_1}=c_1,\qquad \chi(g_2) \pp q{g_2}+\psi(g_2)\pp{q}{g_1}=c_2.
$$
Evaluating the second equation at $g_1=g_2=0$ gives ${c_2}=0$, hence ${c}_1\neq0$ and $q(g_1,g_2) = { c}_1 g_1 +
r(g_2)$ for some smooth function $r(g_2)$.  Then the second condition leads to the following differential equation
for $r$,
\begin{equation}\label{eq:r_diff}
  \chi(g_2) \rp(g_2)=-\psi (g_2) { c}_1.
\end{equation}
But this equation does not admit a smooth solution, unless $\psi(g_2)/\chi(g_2) $ admits a smooth continuation at
$0$. If, for example, $\psi(g_2)=g_2$ and $\chi(g_2)=g_2^2$, then there is no solution $r(g_2)$ to
(\ref{eq:r_diff}), which is smooth in the neighborhood of $0$, hence a system of coordinates in which $\Pi$ takes
the form (\ref{eq:poisson_fic}) does not exist.
\end{example}

% -*- mode: latex; tex-main-file: "art.tex"; -*-
%
\section{Action-angle coordinates for Liouville integrable systems on Poisson manifolds}
\label{sec:aa_coordinates}
In this section we prove the existence of action-angle coordinates in the neighborhood of every standard Liouville
torus of an integrable system on an arbitrary Poisson manifold.

\subsection{Standard Liouville tori of Liouville integrable systems}
We first recall the definition of a Liouville integrable system on a Poisson manifold.
\begin{defn}\label{def:liouville_int}
  Let $(M,\Pi)$ be a Poisson manifold of (maximal) rank $2r$ and of dimension $n$. An $s$-tuplet of functions
  $\F=(f_1,\dots,f_s)$ on $M$ is said to define a \emph{Liouville integrable system} on $(M,\Pi)$ if
  \begin{enumerate}
    \item $f_1,\dots,f_s$ are independent (i.e., their differentials are independent on a dense open subset of $M$);
    \item $f_1,\dots,f_s$ are in involution (pairwise);
    \item $r+s=n$.
  \end{enumerate}
  Viewed as a map, $\F:M\to\R^s$ is called the \emph{momentum map} of $(M,\Pi,\F)$.
\end{defn}
We denote by $M_r$ the open subset of $M$ where the rank of $\Pi$ is equal to~$2r$; points of $M_r$ are called
\emph{regular points} of $M$. We denote by $\U_\F$ the dense open subset of $M$, which consists of all points of
$M$ where the differentials of the elements of $\F$ are linearly independent,
\begin{equation}\label{eq:indep_diff}
  \U_\F:=\set{m\in M\mid \diff_m f_1\we\diff_m f_2\we\dots\we\diff_m f_s\neq0}.
\end{equation}%
On the non-empty open subset $M_r\cap\U_\F$ of $M$ the Hamiltonian vector fields $\X_{f_1},\dots,\X_{f_s}$ define a
distribution $\D$ of rank $r$, since at each point $m$ of $M_r$ the kernel of $\Pi_m$ has dimension $n-2r=s-r$. The
distribution $\D$ is integrable because the vector fields $\X_{f_1},\dots,\X_{f_s}$ pairwise commute,
\begin{equation*}
  \left[\X_{f_i},\X_{f_j}\right]=-\X_{\pb{f_i,f_j}}=0,
\end{equation*}%
for $1\leqs i<j\leqs s$. The integral manifolds of $\D$ are the leaves of a regular foliation, which we denote by
$\FF$; the leaf of $\FF$, passing through $m$, is denoted by $\FF_m$, and is called the \emph{invariant manifold}
of $\F$, through~$m$.  For what follows, we will be uniquely interested in the case in which $\FF_m$ is
compact. According to the classical Liouville theorem, adapted to the case of Poisson manifolds (see \cite[Sect.\
4.3]{adlermoerbekevanhaecke2004} for a proof in the Poisson manifold case), every compact invariant manifold
$\FF_m$ is diffeomorphic to the torus $\T^r:=(\R/\Z)^r$; more precisely, the diffeomorphism can be chosen such that
each of the vector fields $\X_{f_i}$ is sent to a constant (i.e., translation invariant) vector field on
$\T^r$. Such a torus is called a \emph{standard Liouville torus}.

%We consider the subset
%%
%\begin{equation*}
%  \U_\F:=\{m' \in M \mid  \X_{f_1}(m')\we\dots\we\X_{f_r}(m')\neq0\},
%\end{equation*}%
%%
%which is an open subset of $M$, which is contained in $M_r$, and which contains~$m$. For a vector field $\Vf$ and a
%polyvector field $P$ on $M$, we denote the Lie derivative of $P$ with respect to $\Vf$ by $\L_\Vf P$. It is related
%to the Schouten bracket $\LB_S$ by $\L_\Vf P=\lb{\Vf,P}_S$. Since, for $j=1,\dots,r$,
%%
%\begin{equation*}
%  \L_{\X_{f_j}}\(\X_{f_1}\we\dots\we\X_{f_r}\)=\sum_{i=1}^r\X_{f_1}\we\dots\we\left[\X_{f_j},\X_{f_i}\right]
%  \we\dots\we\X_{f_r}=0,
%\end{equation*}%
%%
%$\U_\F$ is invariant for the flow of the Hamiltonian vector fields $\X_{f_1},\dots,\X_{f_r}$.

\subsection{Foliation by standard Liouville tori}
\label{par:fol_by_tori}

As a first step in establishing the existence of action-angle coordinates, we prove that, in some neighborhood of a
standard Liouville torus, the invariant manifolds of an integrable system $(M,\Pi,\F)$ form a trivial torus
fibration.

%\eject

\begin{prop}\label{prp:localdistri}
  Suppose that $\FF_m$ is a standard Liouville torus of an integrable system $(M,\Pi,\F)$ of dimension $n:=\dim M$
  and rank $2r:=\Rk\Pi$. There exists an open subset $U\subset M_r\cap\U_\F$, containing $\FF_m$, and there exists a
  diffeomorphism $\phi:U\simeq \T^r \times B^{n-r}$, which takes the foliation $\FF$ to the foliation, defined by
  the fibers of the canonical projection $p_B:\T^r\times B^{n-r} \to B^{n-r}$, leading to the
  following commutative diagram.
\begin{equation*}
  \begin{diagram}
    \node{\FF_m}\arrow{e,t,J}{}
    \node{U}\arrow{e,t}{\phi}\arrow{e,b}{\simeq}\arrow{s,b}{\F_{\vert U}}
    \node{\T^r\times B^{n-r}}\arrow{sw,r}{p_B}\\
    \node[2]{B^{n-r}}
  \end{diagram}
\end{equation*}%
%
%    \item[$(iv)$] There exists an open subset $U'\subset\U_\F$, containing $\FF_m$, and there exist $n-2r$
%        functions $f_{r+1},\dots,f_{n-r}$ on $U'$, which Poisson commute with each of the functions
%        $f_1,\dots,f_r$, and with the property that the $n-r$ functions $f_1,\dots,f_{n-r}$ have independent
%        differentials at every point of $U'$.
%
\end{prop}
\begin{proof}
%
%We now prove $(ii)\to (iii)$. Thus, we assume that we are given a neighborhood $U'$ of $\FF_m$ in $\U_\F$, and
%Casimir functions $z_1,\dots,z_{n-2r}$ on it, with independent differentials at every point. Put $f_{r+i}:=z_i$ for
%$i=1,\dots,n-2r$. Then, clearly, the $n-r$ functions $f_1,\dots,f_{n-r}$ pairwise Poisson commute on $U'$. To check
%that their differentials are independent, it suffices to check that if $f$ is a Casimir function, then $\diff_mf$
%is linearly independent of $\diff_m f_1,\dots,\diff_m f_r$; but this is clear since $\X_f(m)=0$ while
%$\X_{f_1},\dots,\X_{f_r}$ are linearly independent at $m$. This proves $(ii)\to (iii)$. The proof of $(iii)\to
%(iv)$ is trivial.
%
We first show that the foliation $\FF$, which consists of the maximal integral manifolds of the foliation $\D$,
defined by the integrable vector fields $\X_{f_1},\dots,\X_{f_{s}}$, where $s:=n-r$, coincides with the foliation
$\bar\FF$, defined by the fibers of the submersion
$$
  \bar\F=(f_1,\dots,f_{s}):M_r\cap\U_\F\to\R^{s},
$$
which is the restriction of $\F:M\to\R^{s}$ to $M_r\cap\U_\F$. Since all leaves of $\bar\FF$ and of $\FF$ are
$r$-dimensional, it suffices to show that the two leaves, which pass through an arbitrary point $m\in
M_r\cap\U_\F$, have the same tangent space at $m$. Since $f_1,\dots,f_s$ are pairwise in involution, each of the
vector fields ${\X}_{f_1},\dots,\X_{f_s}$ is tangent to the fibers of $\bar\F$, i.e., to the leaves of
$\bar\FF$. Thus, $T_m\FF\subset T_m\bar\FF$, which implies that both tangent spaces are equal, since they have the
same dimension $r$.

Suppose now that $\FF_m$ is a standard Liouville torus. We show that there exists a neighborhood $U$ of $\FF_m$ and
a diffeomorphism $\phi: U\to\FF_m\times B^{s}$, which sends the foliation $\bar\FF$ $(=\FF)$, restricted to $U$, to
the foliation defined by $p_B$ on $\FF_m\times B^{s}$.  The proof of this fact depends only on the fact that
$\FF_m$ is a compact component of a fiber of a submersion (namely $\bar\F$). Notice that since $\bar\F$ is a
submersion, every point $m'\in\bar\FF_m=\FF_m$ has a neighborhood $U_{m'}$ in~$M$, which is diffeomorphic to the
product of a neighborhood $V_{m'}$ of $m'$ in $\FF_m$ times an open ball $B^{s}_{m'}$, centered at
$\bar\F(m')=\bar\F(m)$ in $\R^{s}$; such a diffeomorphism $\phi_{m'}$, as provided by the implicit function
theorem, is a lifting of $\bar\F$, i.e., it leads to the following commutative diagram:

\begin{equation*}
  \begin{diagram}
    \node{U_{m'}}\arrow{e,t,--}{\phi_{m'}}\arrow{se,b}{\bar\F}\node{V_{m'}\times B_{m'}^{s}}\arrow{s,r}{p_B}\\
    \node[2]{B_{m'}^{s}}
  \end{diagram}
\end{equation*}%
Since $\FF_{m}$ is compact, it is covered by finitely many of the sets $V_{m'}$, say $V_{m_1},\dots,V_{m_\ell}$.
Thus, if every pair of the diffeomorphisms $\phi_{m_1},\dots,\phi_{m_\ell}$ agrees on the intersection of their
domain of definition (whenever non-empty), we can define a global diffeomorphism on a neighborhood $U$ of $\FF_m$,
whose image is the intersection of the concentric balls $B_{m_1}^{s},\dots,B_{m_\ell}^{s}$. In order to ensure that
these diffeomorphisms agree, we need to chose them in a more specific way. This is done by choosing an arbitrary
Riemannian metric on $M$. Using the exponential map, defined by the metric, we can identify a neighborhood of the
zero section in the normal bundle of $\FF_m$, with a neighborhood of $\FF_m$ in $M$; in particular, for every
$m'\in\FF_m$ there exist neighborhoods $U_{m'}$ of $m'$ in $M$ and $V_{m'}$ of $m'$ in $\FF_m$, with smooth maps
$\psi_{m'}:U_{m'}\to V_{m'}$, which have the important virtue that they agree on the intersection of their
domains. Upon shrinking the open subsets $U_{m'}$, if necessary, the maps $\phi_{m'}:=\psi_{m'}\times
(f_{1},\dots,f_{s})$ are a choice of diffeomorphisms, defined on a neighborhood $U$ of $\FF_m$, with the required
properties.
%
%Since $\FF=\bar\FF$, the restriction of both foliations to $U$ have the same leaves, so the result follows.
%
%ir common intersectionin Consider the vector bundle $T\FF_m^{\perp} \to \FF_m$. For any $x \ \FF_m$, denote by $
%p_{x,\epsilon}$ the composition of the exponential map with $\F$ restricted to $E_x \cap B(0,\epsilon)$. For all $x
%\in \FF_m$, there exists a neighborhood $V_x$ of $x \in \FF_m$ and a $\epsilon_x >0$ such that $p_{(y,\epsilon_x)}
%$ is a diffeomorphism onto its image for all $y \in V_x$ and the exponential map, restricted to $\coprod_y
%exp(B(0,\epsilon_x))$, is a diffeomorphism onto its image. By compactness of $\FF_m$, one can cover $\FF_m$ with
%finitely many such open subsets and such an $\epsilon$ can be found globally. Let $T$ be the tubular neighborhood
%of $\FF_m $.  For all $x \in \FF_m$, there exists a neighborhood $W_x$ of $x \in \FF_m$ and a $\eta_x >0$ such that
%$p_{(y,\epsilon_x)} $ contains $ B(0,\eta_x) $ for all $y \in V_x$. By compactness of $\FF_m$, one can cover
%$\FF_m$ with finitely many such open subsets and such an $\eta$ can be found globally. In particular, $
%\F^{-1}\big( B(0,\eta) \big) \subset T$ (at least one of its components).
%
\end{proof}

\begin{cor}\label{cor:casimirs}
  Suppose that $\FF_m$ is a standard Liouville torus of an integrable system $(M,\Pi,\F)$ of dimension $n:=\dim M$
  and rank $2r:=\Rk\Pi$. There exists an open subset $U\subset M_r\cap\U_\F$, containing $\FF_m$, and there exist
  $n-2r$ functions $z_1,\dots,z_{n-2r}$ on $U$ which are Casimir functions of~$\Pi$, and whose differentials are
  independent at every point of $U$.
\end{cor}
\begin{proof}
Let $U\subset M_r\cap\U_\F$ and $\phi$ be as given by proposition \ref{prp:localdistri}. We consider, besides $\D$,
another integrable distribution on $U$: the distribution $\D'$ defined by \emph{all} Hamiltonian vector fields on
$U$; it has rank $2r$ and its leaves are the symplectic leaves of $(U,\Pi)$. Since $\D$ is the distribution,
defined by the Hamiltonian vector fields $\X_{f_1},\dots,\X_{f_s}$, we have that $\D\subset\D'$. Consider the
submersive map $p_B\circ\phi:U\to\T^r\times B^{s}\to B^{s}$, whose fibers are by assumption the leaves of $\FF$,
i.e., the integral manifolds of $\D$ (restricted to~$U$), so that the kernel of $\diff(p_B\circ\phi)$ is precisely
$\D$. The image of $\D'$ by $\diff(p_B\circ\phi)$ is therefore a (smooth) distribution $\D''$ of rank $r$ on
$B^{s}$, which is integrable, since $\D'$ is integrable. The foliation defined by the integral manifolds of $\D''$
is, in the neighborhood of the point $p_B(\phi(m))$, defined by $s-r=n-2r$ independent functions
$z_1',\dots,z'_{n-2r}.$ Pulling them back to $M$, we get functions $z_1,\dots,z_{n-2r}$ on a neighborhood $U$ of
$\FF_m$, with independent differentials on $U$, and they are Casimir functions because they are constant on the
leaves of $\D'$, which are the symplectic leaves of $(U,\Pi)$.
\end{proof}

For Liouville tori in an integrable system, which are not standard, there may not exist a neighborhood on which the
invariant manifolds of the integrable system are locally trivial. We show this in the following example.

\begin{example}
Let $M$ be the product of a M\"obius band with an interval, which is obtained by identifying on $M_0:=\lbrack
{-1},1\rbrack\,\,\times\,\rbrack\!{-1},1\lbrack\,\times\,\R$ in pairs the points $(-1,y,z)$ and $(1,-y,z)$, where
$y$ and $z$ are arbitrary. On $M_0$, consider the vector field $\Vf_0:=\p/\p x$, the Poisson structure
\begin{equation*}
  \Pi_0:=\pp{}x\we\pp{}z,
\end{equation*}
and the function $F:=z$.  The algebra of Casimir functions of $\Pi_0$ consists of all smooth functions on $M_0$
that are independent of $x$ and $z$ (i.e., arbitrary smooth functions in $y$). Clearly, both $\Vf_0$ and $\Pi_0$
and $z$ go down to $M$, yielding a vector field $\Vf$, a Poisson structure $\Pi=\Vf\we\p/\p z$ and a function $z$
on $M$. What does \emph{not} go down to $M$ is the function $y$. In fact, only even functions in $y$ go down and
the algebra of Casimir functions of $\Pi$ is the algebra of even functions in $y$, viewed as functions on $M$. This
remains true if we restrict $M$ to any neighborhood of the central circle $y=z=0$, which is a leaf of the
foliation, defined by the fibers of $F$. Since the differential of an even function in $y$ vanishes at all points
where $y=0$, the central circle is not a standard Liouville torus. Since every neighborhood of the central circle
contains leafs that spin around the M\"obius band twice, the Liouville tori do not form a locally trivial torus
fibration in the neighborhood of the central circle.
\end{example}

\subsection{Standard Liouville tori and Hamiltonian actions}

According to proposition \ref{prp:localdistri}, the study of an integrable system $(M,\Pi,\F)$ in the neighborhood
of a standard Liouville torus amounts to the study of an integrable system $(\T^r\times B^{n-r},\Pi_0,p_B)$, where
$\Pi_0$ is a Poisson structure on $\T^r\times B^{n-r}$ of constant rank $2r$ and the map $p_B:\T^r\times B^{n-r}\to
B^{n-r}$ is the projection on the second factor. We write the latter integrable system in the sequel as
$(\T^r\times B^{s},\Pi,\F)$ and we denote the components of $\F$ by $\F=(f_1,\dots,f_s)$ where $s:=n-r$, as
before. We show in the following lemma that we may assume that the first $r$ vector fields
$\X_{f_1},\dots,\X_{f_r}$ are independent on $\T^r\times B^{s}$, hence span the fibers of $\F$ at each point.
\begin{lemma}\label{lem:indep_torus}
  Let $(\T^r\times B^{s},\Pi,\F)$ be an integrable system, where $\Pi$ has constant rank $2r$ and $\F:\T^r\times
  B^{s}\to B^{s}$ denotes the projection on the second component. Let $m\in \T^r\times \set0$ and suppose that the
  components of $\F=(f_1,\dots,f_s)$ are ordered such that the Hamiltonian vector fields $\X_{f_1},\dots, \X_{f_r}$
  are independent at $m$. There exists a ball $B_0^{s}\subset B^s$, centered at $0$, such that $\X_{f_1},\dots,
  \X_{f_r}$ are independent on $\T^r\times B_0^{s}$.
\end{lemma}
\begin{proof}
We denote by $\L_\V$ the Lie derivative with respect to a vector field $\V$. Since the vector fields $\X_{f_i}$
pairwise commute,
\begin{equation*}
  \L_{\X_{f_j}}\(\X_{f_1}\we\dots\we\X_{f_r}\)=\sum_{i=1}^r\X_{f_1}\we\dots\we\[\X_{f_j},\X_{f_i}\]\we\dots\we\X_{f_r}=0,
\end{equation*}%
for $j=1,\dots,s$. It means that $\X_{f_1}\we\dots\we\X_{f_r}$ is conserved by the flow of each one of the vector
fields $\X_{f_1},\dots,\X_{f_s}$. In particular, if this $r$-vector field is non-vanishing at $m\in\T^r\times\set0$
then it is non-vanishing on the entire integral manifold through $m$ of the distribution $\D$, defined by these
vector fields. Since this integral manifold, which is a torus, is compact, it is actually non-vanishing on a
neighborhood of the integral manifold, which we can choose of the form $\T^r\times B_0^{s}$, where $B_0^{s}\subset
B^s$ is a ball, centered at~$0$.
\end{proof}
%
%nwe have a diffeomorphism $\phi$ between a neighborhood $U$ of the compact leaf $\FF_m$
%of $\FF$ and $\T^r\times B^{n-r}$, and the diffeomorphism respects the foliation, i.e., it sends the standard
%Liouville tori $\FF_{m'}$, for $m'$ in a neighborhood of $m$, to the fibers of $\T^r\times B^{n-r}\to
%B^{n-r}$. Using the diffeomorphism, we may transport the Poisson structure $\Pi$ to $\T^r\times B^{n-r}$, which
%yields a Poisson structure of constant rank $2r$, which we still denote by $\Pi$.
Given an integrable system $(\T^r\times B^{s},\Pi,\F)$, where $\Pi$ has constant rank and $\F=(f_1,\dots,f_s)$ is
the projection on the second component, the Hamiltonian vector fields $\X_{f_i}$ need not be constant on the fibers
of $\F$ (which are tori), and even if they are, they may vary from one fiber to another in the sense that they do
not come from the single action of the torus $\T^r$ on $\T^r\times B^{s}$. We show in the following proposition
how this can be achieved, upon replacing the Hamiltonian vector fields $\X_{f_i}$ by well-chosen linear
combinations, with as coefficients $\F$-basic functions, i.e., functions of the form $\F\circ\l$, where $\l\in
C^\infty(B^s)$; equivalently, smooth functions on $\T^r\times B^{s}$ which are constant on the fibers of $\F$.
\begin{prop}\label{prp:actions1}
  Let $(\T^r\times B^{s},\Pi,\F)$ be an integrable system, where $\Pi$ has constant rank $2r$ and
  $\F=(f_1,\dots,f_s)$ is projection on the second component. Suppose that the $r$ vector fields
  $\X_{f_1},\dots,\X_{f_r}$ are independent at all points of $\T^r\times B^s$. There exists a ball $B_0^s\subset
  B^s$, also centered at $0$, and there exist $\F$-basic functions $\l_i^j\in C^\infty(B_0^{s})$, such that the
  $r$ vector fields $\Y_i:=\sum_{j=1}^r \l^{j}_{i} \X_{f_j}$, $(i=1,\dots,r)$, are the fundamental vector fields of
  a Hamiltonian torus action of $\T^r$ on $\T^r\times B^s_0$.
\end{prop}
The proof uses the following lemma.
\begin{lemma}\label{lem:PoissonCohom}
  Let $\Y$ be a Poisson vector field on a Poisson manifold $(M,\Pi)$ of dimension $n$ and rank $2r$. If $\Y$ is
  tangent to all symplectic leaves of $M$, then $\Y$ is Hamiltonian in the neighborhood of every point $m\in M$
  where the rank of $\Pi$ is $2r$.
\end{lemma}
\begin{proof}
If the rank of $\Pi$ at $m$ is $2r$, so that $m$ is a regular point of $\Pi$, then there exists local coordinates
$(p_1,q_1, \dots, p_r,q_r, z_1, \dots, z_{n-2r}) $ in a neighborhood $U$ of $m$ with respect to which the Poisson
structure $P$ is given by:
$$
  \Pi=\sum_{i=1}^r \pp{}{q_i}\we\pp{}{p_i} .
$$
The vector fields $\pp{}{q_1},\pp{}{p_1},\dots,\pp{}{q_r},\pp{}{p_r}$ span the symplectic leaves of $\Pi$ on
$U$. Therefore, every vector field $\Y$, which is tangent to the symplectic leaves of $\Pi$, is of the form
$$
  \Y =  \sum_{i=1}^r a_i  \pp{}{p_i} + \sum_{i=1}^r  b_i \pp{}{q_i}
$$
for some smooth functions $a_1, \dots, a_r, b_1, \dots b_r$, defined on $U$.  The relation $[\Y,\Pi]=0 $ imposes
the following set of equations to be satisfied for all $i,j=1, \dots, r $:
$$
  \pp{a_i}{q_j} \,=\, \pp{a_j}{q_i}, \, \pp{b_i}{p_j} \,=\, \pp{b_j}{p_i}   \mbox{ and } \pp{a_i}{p_j} \,=\,-\pp{b_i}{q_j}
$$
By the classical Poincar\'e lemma, there exists a function $ h$, defined on $U$, which satisfies, for
$i=1,\dots,r$:
$$
  a_i = -\pp{h}{q_i}  \mbox{ and }  b_i = \pp{h}{p_i}.
$$
Hence,
$$
  \X_h=\sum_{i=1}^r \pp{h}{q_i} \X_{q_i}+\sum_{i=1}^r\pp{h}{p_i}\X_{p_i}+\sum_{k=1}^{n-2r}\pp{h}{z_k}\X_{z_k} = \Y,
$$
which shows that $\Y$ is a Hamiltonian vector field on $U$.
\end{proof}
Now, we can turn our attention to the proof of proposition
\ref{prp:actions1}.
\begin{proof}
The fibers of $\F=(f_1,\dots,f_{s})$ are compact, so for $i=1,\dots,r$, the flow $\Phi^{(i)}_{t_i}$ of the
Hamiltonian vector field $\X_{f_i}$ is complete and we can define a map,
$$
 \begin{array}{ccc}
   {\Phi}:\R^r\times\left(\T^r\times B^{s}\right) &\to &\T^r\times B^{s}\\
   ((t_1,\dots,t_r), m)& \mapsto&\Phi^{(1)}_{t_1}\circ\dots\circ \Phi^{(r)}_{t_r}(m).
 \end{array}
$$
Since the vector fields $\X_{f_i}$ are pairwise commuting, the flows $\Phi^{(i)}_{t_i}$ pairwise commute and $\Phi$
is an action of $\R^r$ on $\T^r\times B^{s}$.  Since the vector fields $\X_{f_1},\dots,\X_{f_r}$ are independent at
all points, the fibers of $\F$, which are $r$-dimensional tori, are the orbits of the action. For $c\in B^{s}$, let
$\Lambda_c$ denote the lattice of $\R^r$, which is the isotropy group of any point in $\F^{-1}(c)$; it is the
period lattice of the action $\Phi$, restricted to $\F^{-1}(c)$. Notice that if $\Lambda_c$ is independent of $c\in
B^{s}$, the action $\Phi$ descends to an action of $\T^r=\R^r/\Lambda_c$ on $\T^r\times B^s$. We will show in Step
1 below that this independence can be assured after applying a diffeomorphism of $\T^r\times B^s_0$ over $B^s_0$,
where $B^s_0$ is a ball, contained in $B^s$, and concentric with it. The proof of this step is essentially the same
as in the symplectic case; it is called \emph{uniformization of the periods}. Steps 2 and 3 below prove
successively that the fundamental vector fields of the obtained torus action are Poisson, respectively Hamiltonian
vector fields.

\step1 The periods of $\Phi$ can be uniformized to obtain a torus action of $\T^r$ on $\T^r\times B^{s}_0$, whose
orbits are the fibers of $\F$ (restricted to $\T^r\times B^{s}_0$).

Let $m_0$ be an arbitrary point of $\F^{-1}(0)$ and choose a basis $(\l_1(0),\dots,\l_r(0))$ for the lattice
$\Lambda_0$. For a fixed $i$, with $1\leqs i\leqs r$, for $m$ in a neighborhood of $m_0$ in $\T^r\times B^{s}$ and
for $L$ in a neighborhood of $\l_i(0)$ in $\R^r$, consider the equation $\Phi(L,m)=m$. Since $\F(\Phi(L,m))=\F(m)$ for
all $L$ and $m$, it is meaningful to write $\Phi(L,m)-m$ and solving the equation $\Phi(L,m)=m$ locally for $L$
amounts to applying the implicit function theorem to the map
\begin{equation*}
  \begin{diagram}
    \node{\R^r\times\left(\T^r\times B^{s}\right)}\arrow[2]{e,t}{{\Phi(L,m)-m}}
    \node[2]{\T^r\times B^{s}}\arrow{e,t}{}
    \node{\T^{r}.}
  \end{diagram}
\end{equation*}%
Since the action is locally free, the Jacobian condition is satisfied and we get by solving for $L$ around
$\l_i(0)$ a smooth $\R^r$-valued function $\l_i(m)$, defined for $m$ in a neighborhood $W_i$ of $m_0$. Doing this
for $i=1,\dots,r$ and setting $W:=\cap_{i=1}^rW_i$, we have that $W$ is a neighborhood of $m_0$, and on $W$ we have
functions $\l_1(m),\dots,\l_r(m)$, with the property that $\Phi(\l_i(m),m)=m$ for all $m\in W$ and for all $1\leqs
i\leqs r$. Thus, $\l_1(m),\dots,\l_r(m)$ belong to the lattice $\Lambda_{\F(m)}$ for all $m\in W$ and they form a
basis when $m=m_0$; by continuity, they form a basis for $\Lambda_{\F(m)}$ for all $m\in W$.

The functions $\l_i$ can be extended to a neighborhood of the torus $\F^{-1}(0)$. In fact, the functions
$\l_i$ are $\F$-basic, hence extend uniquely to $\F$-basic functions on $\F^{-1}(\F(W))$. We will use in the
sequel the same notation $\l_i$ for these extensions and we write $\F^{-1}(\F(W))$ simply as $W$. Using
these functions we define the following smooth map:
\begin{equation}\label{eq:tildephi}
  \begin{array}{lcccc}
   \tilde{\Phi}&:&\R^r\times W &\rightarrow &W \\
   &&((t_1,\dots,t_r), m)& \mapsto & \displaystyle\Phi\left(\sum_{i=1}^r t_i \l_i(m), m\right).
 \end{array}
\end{equation}
Since the functions $\l_i$ are $\F$-basic, the fact that $\Phi$ is an action implies that $\tilde\Phi$ is an
action. The new action has the extra feature that the stabilizer of every point in $W$ is $\Z^r$. Thus,
$\tilde\Phi$ induces an action of $\T^r$ on $W$, which we still denote by $\tilde\Phi$. By shrinking $W$, if
necessary, we may assume that $W$ is of the form $\F^{-1}(B_0^{s})$, where $B_0^{s}$ is an open ball, concentric
with $B^{s}$, and contained in it. Thus we have a torus action
$$
  \begin{array}{lcccc}
   \tilde{\Phi}&:&\T^r\times W &\rightarrow &W \\
   &&((t_1,\dots,t_r), m)& \mapsto & \displaystyle\Phi\left(\sum_{i=1}^r t_i \l_i(m), m\right).
 \end{array}
$$

\smallskip

\step2 The fundamental vector fields of the torus action $\tilde\Phi$ are Poisson vector fields.

We denote by $\Y_1,\dots, \Y_r$ the fundamental vector fields of the torus action~$\tilde{\Phi}$, constructed in step
1. We need to show that $\L_{\Y_i}\Pi=0$, or in terms of the Schouten bracket, that $[\Y_i,\Pi]=0$, for
$i=1,\dots,r$. To do this, we first expand $\Y_i$ in terms of the Hamiltonian vector fields
$\X_{f_1},\dots,\X_{f_r}$: since the action $\tilde{\Phi}$ leaves the fibers of $\F$ invariant and since the
Hamiltonian vector fields $\X_{f_1},\dots,\X_{f_r}$ span the tangent space to these fibers at every point, we can
write
\begin{equation}\label{eq:YX}
  \Y_i=\sum_{j=1}^r \l^{j}_{i} \X_{f_j}.
\end{equation}
Since all Hamiltonian vector fields leave $\Pi$ invariant,
\begin{equation}\label{eq:LiePi}
  \L_{\Y_i}\Pi=[\Y_i,\Pi]=\sum_{j=1}^r\left[\l_i^j\X_{f_j},\Pi\right]= \sum_{j=1}^r\X_{\l_i^j} \we \X_{f_j},
\end{equation}
which we need to show to be equal to zero. Notice that since the coefficients $\l_i$ in the definition of
$\tilde\Phi$ are $\F$-basic, the coefficients $\l^j_i$ are also $\F$-basic, so they are pairwise in involution, and
their Hamiltonian vector fields commute with all Hamiltonian vector fields $\X_{f_k}$. In particular it follows
from~(\ref{eq:LiePi}) that $[\X_{f_k},[\Y_i,\Pi]]=0$ for $k=1,\dots,r$. We derive from it and from~(\ref{eq:YX}),
that $[\Y_i,[\Y_i,\Pi]]=0$, i.e., that the flow of $\Y_i$ preserves $\L_{\Y_i}\Pi$:
\begin{eqnarray} \label{eq:square}
  [\Y_i,[\Y_i,\Pi]]
    &=&\left[\sum_{k=1}^r \l^{k}_{i} \X_{f_k} ,[\Y_i,\Pi]\right]\,=\,
         \sum_{k=1}^r\left[\l^{k}_{i},[\Y_i,\Pi]\right]\we\X_{f_k}\nonumber\\
    &=&\sum_{j,k=1}^r\X_{f_j}\left[\l^{k}_{i}\right]\,\X_{f_k}\we\X_{\l_i^j}+
       \sum_{j,k=1}^r\X_{\l_i^j}\left[\l^{k}_{i}\right]\,\X_{f_k}\we\X_{f_j}\nonumber\\
    &=&0,
\end{eqnarray}
since any two $\F$-basic functions are in involution. Hence, $\L_{\Y_i}^2\Pi =0$. Since $\Y_i$ is a complete vector
field, and has period $1$, we can conclude that $\L_{\Y_i}\Pi=0$ using the following:

\smallskip

\underline{Claim.}$\quad$ If $\Y$ is a complete vector field of period $1$ and $P$ is a bivector field for which
$\L_{\Y}^2 P=0$, then $\L_\Y P=0$.

In order to prove this claim, we let $Q:=\L_\Y P$ and we denote the flow of~$\Y$ by $\Phi_t$. We pick an arbitrary
point $m$ and we show that\footnote{As before, $Q_m$ denotes the bivector $Q$ at the point $m$.}  $Q_m=0$.  We have
for all $t$ that
\begin{equation}\label{eq:dLYP}
  \frac{d}{dt}\left((\Phi_t)_*P_{\Phi_{-t}(m)}\right)=(\Phi_t)_* (\L_\Y P)_{\Phi_{-t}(m)}=(\Phi_t)_* Q_{\Phi_{-t}(m)}=Q_m,
\end{equation}
where we used in the last step that the bivector field $Q$ satisfies $\L_\Y Q=0$.  By integrating (\ref{eq:dLYP}),
$$
  (\Phi_t)_*P_{\Phi_{-t}(m)}=P_m +t Q_m.
$$
Evaluated at $t=1$ this yields $Q_m=0$, since $\Phi_1=\hbox{Id}$, as $\Y$ has period $1$.

\smallskip

\step3 The vector fields $\Y_1,\dots,\Y_r$ are Hamiltonian vector fields
(with respect to commuting Hamiltonian functions).

According to \underline{Step 2}, the vector fields $\Y_1,\dots,\Y_r$ are Poisson vector fields.  Since they are
tangent to the symplectic leaves, according to lemma \ref{lem:PoissonCohom}, there is a neighborhood of $m\in\F_m$
in $W$ that we can assume to be of the form $\Omega_r \times W_s $, with $\Omega_r \subset \T^r , W_s \subset B^s
$, on which the vector fields $\Y_1, \dots,\Y_r$ are Hamiltonian vector fields.  In other words, there exists
functions that we shall denote by $h_1,\dots,h_r $, defined on $W_s$, satisfying the relation $ \Y_i = \X_{h_i}$
for all $i=1, \dots, r$.  It shall be convenient to denote by $W$ again the open subset $\F^{-1}(W_s) $.

Let ${\rm d}\mu $ be a Haar measure on $\T^r $.  For all $m' \in  W$,  we set:
$$
  U_{m'}:=\{t \in \T^r \mid\tilde{\Phi}_t (m') \in W \}
$$
where, for all $t=(t_1, \dots, t_r) \in \T^r $, $ \tilde{\Phi}_t$ is a shorthand for the map $m' \mapsto
\tilde{\Phi}(t_1,\dots, t_r,m')$.  We then define functions $p_i $, $i=1, \dots,r$ on $ W $ by:
$$
  p_i(m'):=\frac{1}{{\rm vol}(U_{m'})} \int_{t \in U_{m'}} h_i\big( \tilde{\Phi}_t (m')  \big) {\rm d} \mu
$$
where ${\rm vol}(U_{m'})$ stands for the volume with respect to the Haar measure.  Their Hamiltonian vector fields
can be computed as follows,:
\begin{eqnarray}
  \X_{p_i}(m')
    &=& \frac{1}{{\rm vol}(U_{m'})}\int_{t \in U_{m'}} \X_{h_i \circ \tilde{\Phi}_t}  (m') {\rm d} \mu  \nonumber\\
    &=& \frac1{{\rm vol}(U_{m'})}\int_{t\in U_{m'}}{\diff}\tilde{\Phi}_t^{-1}\big(\X_{h_i}(\tilde{\Phi}_t(m'))\big){\diff}\mu
    \nonumber\\
    &=&\frac{1}{{\rm vol}(U_{m'})}\int_{t\in U_{m'}}\diff\tilde{\Phi}_t^{-1}\big(\Y_i(\tilde{\Phi}_t(m'))\big){\rm d}\mu
    \nonumber   \\
    &=&\frac{1}{{\rm vol}(U_{m'})} \int_{t \in U_{m'}}  \Y_i (m') {\rm d} \mu \\
    &=& \Y_i(m'),\nonumber  \\
\end{eqnarray}
where the fact that $\Y_i $ is invariant under $\tilde{\Phi}_t $ has been used to go from the third to the fourth
line.  The relation $U_{ \tilde{\Phi}_{t'}(m') } = \tilde{\Phi}_{t'}(U_{m'})$ for all $t' \in \T^r $, and the
invariance property of the Haar measure, imply that the functions $p_1, \dots, p_r $ are invariant under the $\T^r
$-action.  In particular, they are in involution for all $i,j=1,\dots,r$, since
$$
  \{p_i,p_j\} = \Y_j [p_i] =0.
$$
In conclusion, on the open subset $W$, the vector fields $\Y_1,\dots,\Y_r$ are the Hamiltonian vector fields of
the commuting functions $p_1,\dots,p_r $.
\end{proof}

\subsection{The existence of action-angle coordinates}
We are now ready to formulate and prove the action-angle theorem, for standard Liouville tori in Poisson manifolds.
\begin{thm}\label{thm:action-angle}
  Let $(M,\Pi,\F)$ be an integrable system, where $(M,\Pi)$ is a Poisson manifold of dimension $n$ and
  rank~$2r$. Suppose that $\FF_m$ is a standard Liouville torus, where $m\in M_r\cap\U_\F$.  Then there exists
  ${\R} $-valued smooth functions $(p_1,\dots, p_{n-r})$ and $ {\R}/{\Z}$-valued smooth functions
  $({\theta_1},\dots,{\theta_r})$, defined in a neighborhood $U$ of $\FF_m$ such that
  \begin{enumerate}
    \item The functions $(\theta_1,\dots,\theta_r,p_1,\dots,p_{n-r})$ define an isomorphism $U\simeq\T^r\times B^{n-r}$;
    \item The Poisson structure can be written in terms of these coordinates as
      $$
        \Pi=\sum_{i=1}^r\pp{}{\theta_i}\we\pp{}{p_i},
      $$
      in particular the functions $p_{r+1},\dots,p_{n-r}$ are Casimirs of $\Pi$ (restricted to $U$);
    \item The leaves of the surjective submersion $\F=(f_1,\dots,f_{n-r})$ are given by the projection onto the
          second component $\T^r\times B^{n-r}$, in particular, the functions $p_1,\dots,p_{n-r}$ depend on the
          functions $f_1,\dots,f_{n-r}$ only.
  \end{enumerate}
  The functions $\theta_1,\dots,\theta_r$ are called \emph{angle coordinates}, the functions $p_1,$ $\dots,p_r$ are
  called \emph{action coordinates} and the remaining coordinates $p_{r+1},\dots,p_{n-r}$ are called
  \emph{transverse coordinates}.
\end{thm}

\begin{proof}
We denote $s:=n-r$, as before. Since $\FF_m$ is a standard Liouville torus, proposition \ref{prp:localdistri} and
corollary \ref{cor:casimirs} imply that there exist on a neighborhood $U'$ of $\FF_m$ in $M$ on the one hand
Casimir functions $p_{r+1},\dots,p_{s}$ and on the other hand $\F$-basic functions $p_1,\dots,p_r$, such that
$p:=(p_1,\dots,p_{s})$ and $\F$ define the same foliation on $U'$, and such that the Hamiltonian vector fields
$\X_{p_1},\dots,\X_{p_r}$ are the fundamental vector fields of a $\T^r$-action on $U'$, where each of the vector
fields has period 1; the orbits of this torus action are the leaves of the latter foliation. In view of the
Carathéodory-Jacobi-Lie theorem (theorem \ref{thm:localsplitting}), there exist on a neighborhood $U''\subset U'$
of $m$ in $M$, $\R$-valued functions $\theta_1,\dots,\theta_r$ such that
\begin{equation}\label{for:can_loc}
  \Pi = \sum_{j=1}^r\pp{}{\theta_j}\we\pp{}{p_j}.
\end{equation}
On $U''$, $\X_{p_j}=\pp{}{\theta_j}$, for $j=1,\dots,r$; since each of these vector fields has period $1$ on $U'$,
it is natural to view these functions as ${\R}/{\Z}$-valued functions, which we will do without changing the
notation. Notice that the functions $\theta_1, \dots,\theta_r$ are independent and pairwise in involution on $U''$,
as a trivial consequence of (\ref{for:can_loc}). In particular, $\theta_1,\dots,\theta_r,p_1,\dots,p_{s}$ define
local coordinates on $U''$. In these coordinates, the action of $\T^r$ is given by
\begin{equation}\label{for:action_standard}
  (t_1,\dots,t_r)\cdot(\theta_1,\dots,\theta_r,p_1,\dots,p_{s})=(\theta_1+t_1,\dots,\theta_r+t_r,p_1,\dots,p_{s}),
\end{equation}
so that the functions $\theta_i$ uniquely extend to smooth $\R/\Z$-valued functions satisfying
(\ref{for:action_standard}), on $U:=\F^{-1}(\F(U''))$, which is an open subset of $\F_m$ in $M$; the extended
functions are still denoted by $\theta_i$. It is clear that $\pb{\theta_i,p_j}=\delta_i^j$ on $U$, for all $i,j=
1,\dots,r$.  Combined with the Jacobi identity, this leads to
$$
  \X_{p_k}[\pb{\theta_i,\theta_j}]=\pb{\pb{\theta_i,\theta_j},p_k}=\pb{\theta_i,\delta_j^k}-\pb{\theta_j,\delta_i^k}=0,
$$
which shows that the Poisson brackets $\pb{\theta_i,\theta_j} $ are invariant under the $\T$-action; but the latter
vanish on $U''$, hence these brackets vanish on all of $U$, and we may conclude that on $U$, the functions
$(\theta_1,\dots,\theta_r,p_1,\dots,p_{s})$ have independent differentials, so they define a diffeomorphism to
$\T^r\times B^{s}$ where $B^{s}$ is a (small) ball with center $0$, and that the Poisson structure takes in
terms of these coordinates the canonical form (\ref{for:can_loc}), as required.
\end{proof}

The results of the present section can be applied in particular for
a well-known integrable system constructed on a regular coadjoint
orbit ${\mathcal O} $ of ${\mathfrak u}(n)^* $, namely the
Gelfand-Cetlin integrable system, for which action-angle coordinates
are computed explicitly in \cite{guilleminandsternberg} and
\cite{Giacobbe}.  This system can be seen in the Poisson setting, as
follows.  Dualizing the increasing sequence of Lie algebra
inclusions:
$$
  {\mathfrak u}(1) \subset \dots\subset  {\mathfrak u}(n-1) \subset {\mathfrak u}(n)
$$
(where ${\mathfrak u}(k) $ is considered as the left-upper diagonal block of ${\mathfrak u}(k+1) $ for $k=1,\dots,
n-1$), we get a sequence of surjective Poisson maps:
\begin{equation*}
  \begin{diagram}
    \node{{\mathfrak u}(n)^*}\arrow{e,A}
    \node{{\mathfrak u}(n-1)^*}\arrow{e,A}
    \node{\cdots}\arrow{e,A}
    \node{{\mathfrak u}(1)^*}
  \end{diagram}
\end{equation*}%
%
%
%$$
%  {\mathfrak u}(n)^* \to {\mathfrak u}(n-1)^* \to \dots \mapsto  {\mathfrak u}(1)^*.
%$$
%
The family of functions on ${\mathfrak u}(n)^* $ obtained by pulling-back generators of the Casimir algebras of all
the ${\mathfrak u}(k)^* $ for $k=1, \dots, n$ yields a Liouville integrable system on ${\mathfrak u}(n)^*$.  For
particular generators, its restriction to an open subset of ${\mathcal O} $ gives the Gelfand-Cetlin system.  The
invariant manifold is compact, so that theorem \ref{thm:action-angle} can be applied and gives the existence of
action-angle coordinates, defined not only in a neighborhood of the invariant manifold in ${\mathcal O} $, but in a
neighborhood of the invariant manifold in the ambient space~${\mathfrak u}(n) $. The restriction of these action
and angle coordinates to one symplectic leaf~$\mathcal O$ will give action-angle coordinates on $\mathcal O$, as in
\cite{guilleminandsternberg} or \cite{Giacobbe}.

% -*- mode: latex; tex-main-file: "art.tex"; -*-
%
\section{Action-angle coordinates for non-commutative integrable systems on Poisson manifolds}
  \label{sec:non-com}
In this section, we prove the existence of action-angle coordinates in a neighborhood of a compact invariant
manifold in the very general context of non-commutative integrable systems.
\subsection{Non-commutative integrable systems}\label{par:non-com}
We first define precisely what we mean by a non-commutative integrable system on a Poisson manifold, since the
definitions in the literature \cite{bolsinov,AF, FA} are only given in the case of a symplectic manifold. See the
appendix for a more intrinsic version of this definition.
\begin{defn}\label{def:non-com}
Let $(M,\Pi)$ be a Poisson manifold. An $s$-tuple of functions
$\F=(f_1,\dots,f_s)$ is said to be a \emph{non-commutative
integrable system} of rank $r$ on $(M,\Pi)$ if
\begin{enumerate}
  \item[(1)] $f_1,\dots,f_s$ are independent (i.e. their differentials are independent on a dense open subset of $M$);
  \item[(2)] The functions $f_1,\dots,f_r$ are in involution with the functions $f_1,\dots,f_s$;
  \item[(3)] $r+s =\dim M$;
  \item[(4)] The Hamiltonian vector fields of the functions $f_1, \dots,f_r$ are linearly independent at some point of $M$.
\end{enumerate}
We denote the subset of $M$ where the differentials $\diff f_1,\dots,\diff f_s$ (resp.\ where the Hamiltonian
vector fields $\X_{f_1}, \dots,\X_{f_r}$) are independent by $\U_\F$ (resp.\ by $M_{\F, r}$). Notice that
$2r\leqs\Rk\Pi$, as a consequence of (4).
\end{defn}
If $(M,\Pi,\F)$ is a Liouville integrable system (definition
\ref{def:liouville_int}), then it is clear that the components
$(f_1,\dots,f_s)$ of $\F$ can be ordered such that $\F$ is a
non-commutative integrable system of rank $\frac12\Rk\Pi$. Thus, the
notion of a non-commutative integrable system on a Poisson manifold
$(M,\Pi)$ generalizes the notion of a Liouville integrable system on
$(M,\Pi)$. For simplicity, we often refer in this section to the
case of a Liouville integrable system as the \emph{commutative
case}.

\subsection{Standard Liouville tori for non-commutative integrable systems}
Let $\F$ be a non-commutative integrable system of rank $r$ on a Poisson manifold $(M,\Pi)$ of dimension $n$.  The
open subsets $\U_\F$ and $M_{\F,r}$ are preserved by the flow of each of the vector fields
$\X_{f_1},\dots,\X_{f_r}$ since each of the functions $f_1,\dots,f_r$ is in involution with all the functions
$f_1,\dots,f_s$. On the non-empty open subset $M_{\F,r}\cap\U_\F $ of $M$, the Hamiltonian vector fields
$\X_{f_1},\dots,\X_{f_r}$ define a (regular) distribution $\D$ of rank $r$.  Since the vector fields
$\X_{f_1},\dots,\X_{f_r}$ commute pairwise, the distribution $\D$ is integrable, and its integral manifolds are the
leaves of a (regular) foliation ${\mathcal F} $. The leaf through $m \in M$ is denoted by $\FF_m $, and called the
\emph{invariant manifold through $m$} of $\F$. As in the commutative case, we are only interested in the case where
$\FF_m$ is compact.  Under this assumption, $\FF_m$ is a compact $r$-dimensional manifold, equipped with $r$
independent commuting vector fields, hence it is diffeomorphic to an $r$-dimensional torus $\T^r$; then $\FF_m$ is
called a \emph{standard Liouville torus} of $\F $. Proposition \ref{prp:localdistri} takes in the general situation
of a non-commutative integrable system formally the same form, but with the understanding that $r$ now stands for
the rank of $\F$ (rather than half the rank of the Poisson structure), as stated in the following
proposition\footnote{Recall that $B^{n-r}$ is a ball of dimension $n-r$.}.
\begin{prop}\label{prp:localdistri'}
  Suppose that $\FF_m$ is a standard Liouville torus of a non-commutative integrable system $\F$ of rank $r$ on an
  $n$-dimensional Poisson manifold $(M,\Pi)$. There exists an open subset $U\subset M_{\F,r}\cap\U_\F$, containing
  $\FF_m$, and there exists a diffeomorphism $\phi:U\simeq \T^r \times B^{n-r}$, which takes the foliation $\FF$ to
  the foliation, defined by the fibers of the canonical projection $p_B:\T^r\times B^{n-r} \to B^{n-r}$, leading to
  the following commutative diagram.
\begin{equation*}
  \begin{diagram}
    \node{\FF_m}\arrow{e,t,J}{}
    \node{U}\arrow{e,t}{\phi}\arrow{e,b}{\simeq}\arrow{s,b}{\F_{\vert U}}
    \node{\T^r\times B^{n-r}}\arrow{sw,r}{p_B}\\
    \node[2]{B^{n-r}}
  \end{diagram}
\end{equation*}%
\end{prop}
\subsection{Standard Liouville tori and Hamiltonian actions}
According to proposition \ref{prp:localdistri'}, the study of a non-commutative integrable system $(M,\Pi,\F)$ of
rank $r$ in the neighborhood of a standard Liouville torus amounts to the study of the non-commutative integrable
system $(\T^r\times B^{n-r},\Pi_0,p_B)$ of rank $r$, where $\Pi_0$ is a Poisson structure on $\T^r\times B^{n-r}$
and the map $p_B:\T^r \times B^{n-r}\to B^{n-r}$ is the projection onto the second factor. We write the latter
integrable system in the sequel as $(\T^r\times B^{s},\Pi,\F)$ and we denote the components of $\F$ by
$\F=(f_1,\dots,f_s)$ where $s:=n-r$, as before. We may assume that the first $r$ vector fields
$\X_{f_1},\dots,\X_{f_r}$ are independent on $\T^r\times B^{s}$, as shown in the following lemma, the proof of
which goes along the same lines as the proof of lemma \ref{lem:indep_torus}.
\begin{lemma}
  Let $(\T^r\times B^{s},\Pi,\F)$ be a non-commutative integrable system of rank $r$, where $\F:\T^r\times B^{s}\to
  B^{s}$ denotes the projection onto the second component. Let $m\in \T^r\times \set0$ and suppose that the
  Hamiltonian vector fields $\X_{f_1},\dots, \X_{f_r}$ are independent at $m$. There exists a ball $B_0^{s}\subset
  B^s$, centered at~$0$, such that $\X_{f_1},\dots, \X_{f_r}$ are linearly independent at every point of
  $\T^r\times B_0^{s}$.
\end{lemma}
One useful consequence of the fact that the Hamiltonian vector fields $\X_{f_1},\dots,\X_{f_r}$ are independent on
$M:=\T^r\times B^{s}$ is that a function $g\in C^\infty(M)$ is $\F$-basic if and only if $\X_{f_i}[g]=0$ for
$i=1,\dots,r$. Indeed, $g$ is $\F$-basic if and only $g$ is constant on all fibers of $\F$, and all tangent spaces
to these fibers are spanned by the vector fields $\X_{f_1},\dots,\X_{f_r}$.

\smallskip

We now come to an important difference between the commutative and the non-commutative case, which is related to
the nature of the map $\F$. In the commutative case, two $\F$-basic functions on $\T^r\times B^{s}$ are in
involution, $\pb{g\circ\F,h\circ\F}=0$ for all $g,h\in C^\infty(B^s)$. Said differently,
$$
  \F:(\T^r\times B^{s},\Pi)\to (B^s,0),
$$
is a Poisson map, where $B^s$ is equipped with the trivial Poisson structure. The generalization to the
non-commutative case is that $B^s$ admits a Poisson structure (non-zero in general), such that $\F$ is a Poisson
map. This Poisson structure is constructed by the following (classical) trick: for every pair of functions $g,h\in
C^\infty(B^s)$ we have in view of the Jacobi identify, for all $i =1,\dots, r$,
$$
  \X_{f_i}[\pb{g\circ\F,h\circ\F}]=\pb{\X_{f_i}[g\circ\F],h\circ\F}+\pb{g\circ\F,\X_{f_i}[h\circ\F]}=0,
$$
so that $\pb{g \circ \F, h \circ \F}$ is $\F$-basic, namely $\pb{g\circ\F,h\circ\F}=\pb{g,h}_B\circ\F$ for some
function $\pb{g,h}_B \in C^\infty(B^{s})$. It is clear that this defines a Poisson structure $\Pi_B=\PB_B$ on $B^s$
and that
$$
  \F:(\T^r\times B^{s},\Pi)\to (B^s,\Pi_B)
$$
is a Poisson map. This Poisson structure leads to a special class of $\F$-basic functions, which play an important
role in the non-commutative case, defined as follows.
\begin{defn}
  A smooth function $h$ on $\T^r\times B^{s}$ is said to be a \emph{Casimir-basic function}, or simply a
  \emph{Cas-basic function} if there exists a Casimir function $g$ on $(B^s,\Pi_B)$, such that $h=g\circ\F$.
\end{defn}
A characterization and the main properties of Cas-basic functions are given in the following proposition.
\begin{prop}\label{Casimir-basic}
  Let $\F$ be a non-commutative integrable system on a Poisson manifold $(M,\Pi)$, where $M=\T^r\times B^{s}$ and
  $\F$ is projection on the second component. It is assumed that the Hamiltonian vector fields
  $\X_{f_1},\dots,\X_{f_r}$ are independent at every point of $M$.
  \begin{enumerate}
    \item If $g\in C^\infty(M)$, then $g$ is Cas-basic if and only $g$ is in involution with every function which is
           constant on the fibers of~$\F$;
    \item Every pair of Cas-basic functions on $M$ is in involution;
    \item If $g$ is Cas-basic, then its Hamiltonian vector field $\X_g$ on $M$ is of the form
          $\X_F=\sum_{i=1}^r\psi_i\X_{f_i}$, where each $\psi_i$ is a Cas-basic function on $M$.
  \end{enumerate}
\end{prop}
\begin{proof}
Suppose that $g\in C^\infty(M)$ is in involution with every function which is constant on all fibers of $\F$. Then
$\X_{f_i}[g]=\pb{g,f_i}=0$ for $i=1,\dots,r$, hence $g$ is $\F$-basic, $g=h\circ \F$ for some function $h$ on
$B^s$. If $k\in C^\infty(B^s)$, then $k\circ\F$ is constant on the fibers of $\F$, so that
\begin{equation*}
  \pb{h,k}_B\circ\F=\pb{g,k\circ\F}=0,
\end{equation*}%
where we have used that $\F$ is a Poisson map. It follows that $\pb{h,k}_B=0$ for all functions $k$ on $B^s$, hence
that $g$ ($=h\circ\F$) is Cas-basic. This shows one implication of (1), the other one is clear. (2) is an easy
consequence of~(1). Consider now the Hamiltonian vector field $\X_g$ of a Cas-basic function $g$ on~$M$. In view of
(1), $\X_g[h]=\pb{h,g}=0$ for every function $h$ which is constant on the fibers of $\F$, hence $\X_g$ is tangent
to the fibers of $\F$. Since the fibers of $\F$ are spanned at every point by the Hamiltonian vector fields
$\X_{f_1},\dots,\X_{f_r}$, there exist smooth functions $\psi_1,\dots,\psi_r$ on $M$, such that
\begin{equation*}
  \X_g=\sum_{i=1}^r\psi_i\X_{f_i}.
\end{equation*}%
The functions $\psi_i$ are $\F$-basic, because $\X_{h}[\psi_i]=0$ for every function $h$ which is constant on the
fibers of $\F$. Indeed, for such a function $h$, we have that $\pb{g,h}=0$ and $\pb{f_i,h}=0$ for $i=1,\dots,r$, so
that
\begin{equation*}
  0=\X_{\pb{g,h}}=[\X_h,\X_g]=\sum_{i=1}^r\(\X_{h}[\psi_i]\X_{f_i}+\psi_i\lb{\X_h,\X_{f_i}}\)=\sum_{i=1}^r\X_{h}[\psi_i]\X_{f_i}
\end{equation*}%
and the result follows from the independence of $\X_{f_1},\dots,\X_{f_r}$.
\end{proof}
Now, we can give a proposition that generalizes proposition
\ref{prp:localdistri} to the non-commutative setting, which has
formally the same shape up to the fact that $r$, formerly half of
the rank of the Poisson structure $ \Pi$, stands now for rank of the
non-commutative integrable system, and up to the fact that the
functions $\lambda_i^j $ that appear below, are now proved to be
Cas-basic, and not simply $\F$-basic.
\begin{prop}\label{prp:actions2}
Let $(\T^r \times B^{s},\Pi,\F)$ be an non-commutative integrable system of rank $r$, where $\F=(f_1,\dots,f_s)$ is
projection on the second component. There exists a ball $B_0^s\subset B^s$, also centered at $0$, and there exist
Cas-basic functions $\l_i^j$, such that the $r$ vector fields $\Y_i:=\sum_{j=1}^r \l^{j}_{i} \X_{f_j}$,
$(i=1,\dots,r)$, are the fundamental vector fields of a Hamiltonian torus action of $\T^r$ on $\T^r\times B^s_0$.
\end{prop}
We can now turn our attention to the proof of proposition
\ref{prp:actions2}.
\begin{proof}
As in \underline{Step 1} of the proof of proposition \ref{prp:actions1},we obtain the existence of a family of
$\R^r$-valued $\F$-basic functions $\lambda_1,\dots, \lambda_r $ such that $\tilde{\Phi}$, defined as in
(\ref{eq:tildephi}), induces a $\T^r$-action on $\T^r \times B^s_0$, where $B_0^s$ is an $s$-dimensional ball,
contained in $B^s$. As in \underline{Step 2}, we expand the fundamental vector fields $\Y_1,\dots,\Y_r$ of the
action in terms of the Hamiltonian vector fields $\X_{f_1},\dots,\X_{f_r}$,
\begin{equation*}
    \Y_i=\sum_{j=1}^r \l^{j}_{i} \X_{f_j}.
\end{equation*}%
The proof that the vector fields $\Y_i$ are Poisson vector fields requires an extra argument: we show that the
relations $\X_{f_i} [\lambda_i^k] =0=\X_{\lambda_i^j} [\lambda_i^k]=0$ which were used in (\ref{eq:square}) still
hold, by showing that the functions $\l_i^j$ are Cas-basic (recall that Cas-basic functions are in involution). To
do this, it suffices to show that if a vector field on $\T^r\times B^s_0$ of the form
$Z=\sum_{i=1}^r\psi_i\X_{f_i}$ is periodic of period $1$, then each of the coefficients $\psi_i$ is Cas-basic. Let
$Z$ be such a vector field and consider
\begin{equation*}
  Z_0:=\sum_{i=1}^r\psi_i(m)\X_{f_i},
\end{equation*}%
where $m$ is an arbitrary point in $\T^r\times B^s_0$. Then the restriction of $Z_0$ to $\F^{-1}(\F(m))$ is
periodic of period $1$. Let $h$ be an $\F$-basic function on $\T^r\times B^s_0$, and let us denote the (local) flow
of $\X_h$ by $\Phi_t$. Since
\begin{equation*}
  \lb{\X_h,Z_0}=\sum_{i=1}^r\psi_i(m)\lb{\X_h,\X_{f_i}}=0,
\end{equation*}%
for $\vert t\vert$ sufficiently small, the flow of $Z_0$ starting from $\Phi_t(m)$ is also periodic of period $1$.
Since the coefficients of $Z$ are the unique continuous functions such that $Z=Z_0$ on $\F^{-1}(\F(m))$ and such
that the flow of $Z$ from every point has period $1$, it follows that $\psi_i(\Phi_t(m))=\psi_i(m)$ for $\vert
t\vert$ sufficiently small. Taking the limit $t\mapsto 0$ yields that $\X_h[\psi_i]=0$ for every $\F$-basic
function on~$\T^r\times B^s_0$. Thus, $\psi_i$ is Cas-basic, for $i=1,\dots,r$.
so that the proof of \underline{Step 2} remains valid,
amounting to the fact that the vector fields $\Y_i$ on $W$ are Poisson vector fields, $\L_{\Y_i}\Pi=0$, which leads in
view of (\ref{eq:LiePi}) to
\begin{equation}\label{eqn:thispoisson}
  \sum_{j=1}^r \X_{\l^{j}_{i}}\we\X_{f_j}=0.
\end{equation}
We show that these vector fields are Hamiltonian, where the Hamiltonians can be taken as a $\F$-basic functions.
The key point is that all coordinates which appear all along this step should now be taken with respect to coordinates adapted to $f_1,
\dots, f_r$.  More precisely, we choose some $ m \in \T^r \times B^s_0$ in $\F^{-1} (0)$, and we construct in some
neighborhood $W_0' $ of $m$ a system of coordinates
$$
  (f_1,g_1, \dots, f_{r}, g_r, z_1, \dots, z_{n-2r})
$$
in which the Poisson structure takes the form given in equation (\ref{eq:thm_split}). Of course, the functions
$z_1,\dots,z_{n-2r}$ are $\F$-basic again (and therefore depend on $f_1, \dots, f_s $ only), so that they can be
defined in $p^{-1}(p(W_0'))$, an open subset which we also call $W_0'$ for the sake of simplicity.  As before, we
make no notational distinction between the functions $f_1, \dots, f_r,z_1, \dots, z_{n-2r} $, considered as
functions on $p(W_0') \subset B^{n-r} $, and the functions $f_1, \dots, f_r, z_1, \dots, z_{n -2r} $ themselves,
defined on $ W_0'$.

In view of Proposition \ref{Casimir-basic}(3), in the previous system of coordinates, we have that,
since the functions $ \lambda_{i}^j$ are Cas-basic,
$$
  \X_{\lambda_{i}^j} = \sum_{k=1}^r \big( \pp{\mu_i^j}{{f}_k} \circ \F \big ) \, \X_{f_k},
$$
Hence,  (\ref{eqn:thispoisson}) gives
%To do this, we write (\ref{eqn:thispoisson}) in terms of $\F$-basic functions, which we do by writing the $\F$-basic
%functions $\l_i^j$ on $W$, explicily in terms of functions on $\F(W)$: there exist smooth functions $\mu_i^j$ on
%$\F(W)$, such that $\l^j_i=\mu^j_i\circ \F$, which upon substitution in (\ref{eqn:thispoisson}) yields
%
\begin{equation} \label{eq:partial_deriv}
  \sum_{1\leqs j<k\leqs r}\left(\left(\pp{\mu^j_i}{f_k}-\pp{\mu^k_i}{f_j}\right)\circ \F\right)\,\X_{f_j}\we\X_{f_k}=0.
\end{equation}
where $\mu_i^j$ is defined by $\l^j_i=\mu^j_i\circ \F$.
Since the vectors fields $\X_{f_1},\dots,\X_{f_r}$ are linearly independent at all points of $W$, all coefficients
above vanish and we get, for every $i,j,k \in\{1,\dots,r\}$:
\begin{equation}\label{eqn:Poincare}
  \pp{\mu^j_i}{f_k}=\pp{\mu^k_i}{f_j}.
\end{equation}
As in the proof of the classical Poincar\'{e} lemma, the functions $b_1,\dots,b_r$ on $\F(W)$, defined by
\begin{equation}\label{eq:bi}
  b_i=b_i(f_1,\dots,f_{s}):=\sum_{j=1}^r\int_{t=0}^1\mu^j_i(tf_1,\dots,tf_r,z_{1},\dots,z_{s-r}) f_j
\end{equation}
satisfy
\begin{equation}\label{eq:poincare}
  \mu^j_i=\pp{b_i}{f_j},
\end{equation}
for all $1\leqs i,j\leqs r$. It leads to the $\F$-basic functions $p_1,\dots,p_r$, defined by
$p_i:=b_i\circ\F$, for $i=1,\dots,r$.

Since the Poisson structure on $B^s$ depends only on the variables $ z_1, \dots, z_s$,
  the map $(f_1, \dots, f_r,z_1, \dots, z_{s-r}) \mapsto (t f_1, \dots, t f_r,z_1, \dots, z_{s-r}) $
is Poisson for all  $t\in [0,1]$. Therefore, the function $ \mu^j_i(tf_1,\dots,tf_r,z_{1},\dots,z_{s-r})$ is a Casimir function,
and the functions $b_1, \dots, b_r $, which are obtained by integration (w.r.t.\ $t$) of these
functions,  are also Casimir functions.
Hence, the functions defined by $p_i:=b_i\circ\F$ are Cas-basic. The Hamiltonian vector field of $p_i$ is, in view of
in view of Proposition \ref{Casimir-basic}(3), (\ref{eq:poincare}) and (\ref{eq:YX}), given by
\begin{equation}\label{eq:rhoi}
  \X_{p_i}=\X_{b_i\circ \F}=\sum_{j=1}^r\left(\pp{b_i}{f_j}\circ\F\right)\,\X_{f_j}
              =\sum_{j=1}^r\left(\mu_i^j\circ \F\right)\,\X_{f_j}=\sum_{j=1}^r\l_i^j\X_{f_j}=\Y_i.
\end{equation}%
This shows that each one of the vector fields $\Y_i$ is a Hamiltonian vector field on $W$.
This completes the proof.
%
%so that (\ref{eq:partial_deriv}) holds again, and the functions $b_1,\dots,b_r$, defined as in (\ref{eq:bi}),
%satisfy the relation (\ref{eq:poincare}). Since these functions are obtained by integrating (w.r.t.\ $t$) the
%Casimir functions $\mu_i^j(t{f}_1,\dots,t{f}_r,{z}_{1},\dots,{z}_{n-2r})$, the functions $p_1,\dots, p_r$, defined
%by $p_i:=b_i\circ\F$ are Cas-basic, so that the relation (\ref{eq:rhoi}) is still valid. This completes the proof.
\end{proof}

\subsection{The existence of action-angle coordinates}
We finally get to the action-angle theorem,  for standard Liouville tori of a non-commutative integrable system.
\begin{thm}\label{thm:action-angle2}
  Let $(M,\Pi)$ be a Poisson manifold of dimension $n$, equipped with a non-commutative integrable system
  $\F=(f_1,\dots,f_s)$ of rank $r$, and suppose that $\FF_m$ is a standard Liouville torus, where $m\in
  M_{\F,r}\cap\U_\F$.  Then there exist ${\R} $-valued smooth functions $(p_1,\dots, p_r,z_1, \dots, z_{s-r} )$ and
  $ {\R}/{\Z}$-valued smooth functions $({\theta_1},\dots,{\theta_r})$, defined in a neighborhood $U$ of $\FF_m$,
  and functions such that
  \begin{enumerate}
    \item The functions $(\theta_1,\dots,\theta_r,p_1,\dots,p_{r},z_1 ,\dots, z_{s-r} )$ define an isomorphism
            $U\simeq\T^r\times B^{s}$;
    \item The Poisson structure can be written in terms of these coordinates as
      $$
        \Pi=\sum_{i=1}^r \pp{}{\theta_i}\we\pp{}{p_i} + \sum_{k,l=1}^{s-r} \phi_{k,l}(z) \pp{}{z_k}\we\pp{}{z_l} ;
      $$
    \item The leaves of the surjective submersion $\F=(f_1,\dots,f_{s})$ are given by the projection onto the
          second component $\T^r \times B^{s}$, in particular, the functions $p_1,\dots,p_r,z_1, \dots, z_{s-r} $
          depend on the functions $f_1,\dots,f_s$ only.
  \end{enumerate}
  The functions $\theta_1,\dots,\theta_{r}$ are called \emph{angle coordinates}, the functions $p_1,\dots,p_r$
  are called \emph{action coordinates} and the remaining coordinates $z_1,\dots,z_{s-r}$ are called
  \emph{transverse coordinates}.
\end{thm}

\begin{proof}
Conditions (1) and (2), in view of lemma \ref{lem:indep_torus}, propositions \ref{prp:localdistri'} and
\ref{prp:actions1} imply that there exist on a neighborhood $U'$ of $\FF_m$ in $M$ on the one hand $\F$-basic
functions $z_{1},\dots,z_{s-r}$ and on the other hand Cas-basic functions $p_1,\dots,p_r$, such that
  $$
    p:=(p_1,\dots,p_{r},z_{1},\dots,z_{s-r})
  $$
and $\F$ define the same foliation on $U'$, and such that the Hamiltonian vector fields $\X_{p_1},\dots,\X_{p_r}$
are the fundamental vector fields of a $\T^r$-action on $U'$, where each has period 1; the orbits of this torus
action are the leaves of the latter foliation. In view of theorem \ref{thm:localsplitting}, there exist on a
neighborhood $U''\subset U'$ of $m$ in $M$, $\R$-valued functions $\theta_1,\dots,\theta_r$ such that
\begin{equation*}
  \Pi = \sum_{j=1}^r\pp{}{\theta_j}\we\pp{}{p_j} + \sum_{k,l=1}^{s-r} \phi_{k,l}(z)  \pp{}{z_k}\we\pp{}{z_l}.
\end{equation*}%
The end of the proof goes along the same lines as the end of the proof of theorem \ref{thm:action-angle}.
\end{proof}

%
% -*- mode: latex; tex-main-file: "art.tex"; -*-
%
\section{Appendix: non-commutative integrability on Poisson manifolds}
In the symplectic context, the terms superintegrability, non-commutative integrability, degenerate integrability,
generalized Liouville integrability and Mischenko-Fomenko integrability refer to the case when the Hamiltonian flow
admits more independent constants of motions than half the dimension of the symplectic manifold
\cite{bolsinov,FA,AF,Nekhroshev,MF}. All these names correspond to notions which are equivalent, at least
locally. Similarly, the definition of a non-commutative integrable system on a Poisson manifold, which we have
given in section \ref{sec:non-com} (definition \ref{def:non-com}) admits different locally equivalent formulations,
which each have their own flavor. We illustrate this in this appendix, by giving an abstract geometrical
formulation in terms of foliations, and a concrete geometrical formulation in terms of Poisson maps.

For both geometrical formulations, the notion of polarity in Poisson geometry is a key element. Let $m$ be an
arbitrary point of a Poisson manifold $(M,\Pi)$.  The \emph{polar} of a subspace $\Sigma\subset T^*_mM$ is the
subspace $\Sigma^{pol}\subset T^*_m M$, defined by
$$
   \Sigma^{pol}:= \{\xi \in T^*_mM \mid \Pi_m (\xi, \Sigma) =0 \} .
$$
Notice that the polar of $\Sigma^{pol}$ can be strictly larger than
$\Sigma$, because $\Pi_m$ may have a non-trivial kernel. When
$\Sigma=\Sigma^{pol}$, we say that $\Sigma$ is a Lagrangian
subspace.

Let $\FF$ and $\GG$ be two foliations on the same Poisson manifold $(M,\Pi)$. For $m\in M$ we denote by
$T_m^\perp\FF$ the subspace of $T_m^*M$, consisting of all covectors which annihilate $T_m\FF$, the tangent space
to the leaf of $\FF$, passing through~$M$. If $\FF$ is defined around $m$ by functions $f_1,\dots,f_s$, then
$T_m^\perp\FF$ is spanned by $\diff_m f_1,\dots,\diff_mf_s$. We say that $\FF$ is \emph{polar} to $\GG$ if
$T^\perp_m \FF =(T^\perp_m\GG)^{pol}$, for every $m \in M$; also, $\FF$ is said to be a \emph{Lagrangian foliation}
if $\FF$ is polar to $\FF$.

\begin{defn}\label{def:abstract-non-Com}
Let $(M,\Pi) $ be a Poisson manifold.  An \emph{abstract non-commutative integrable system} of rank $r$ is given by
a pair $(\FF,\GG) $ of foliations on $M$, satisfying
\begin{enumerate}
  \item $\FF$ is of rank $r$ and $\GG$ is of corank $r$;
  \item $\FF$ is contained in $\GG $ (i.e., each leaf of $\FF$ is contained in a leaf of $\GG$);
  \item $\FF$ is polar to $\GG$.
\end{enumerate}
\end{defn}

This definition generalizes the definition of an \emph{abstract} integrable system on $(M,\Pi)$, which is simply a
Lagrangian foliation $\FF$ on $M$.

\smallskip

In the following proposition we prove the precise relation between definitions \ref{def:non-com} and
\ref{def:abstract-non-Com}.

\begin{prop}
  Let $(M,\Pi) $ be a Poisson manifold.
  \begin{enumerate}
    \item If $\F = (f_1, \dots, f_{s}) $ is a non-commutative integrable system of rank $r$ on $(M,\Pi)$, then on
           $\U_\F \cap M_{\F,r}$ the pair of foliations $(\FF,\GG)$, defined by $\FF :={\rm fol}(f_1,\dots,f_s)$ and
           $\GG:={\rm fol}(f_1,\dots,f_r)$ is an abstract non-commutative integrable system of rank $r$;
    \item Given $(\FF,\GG)$ an abstract non-commutative integrable system of rank~$r$ on $(M,\Pi)$, there exists
          for every $m$ in~$M$ a neighborhood $ U $ in $M$, and functions $\F=(f_1,\dots,f_s)$ on $U$, such that
          $\F$ is a non-commutative integrable system of rank $r$ on $U$.
  \end{enumerate}
\end{prop}
\begin{proof}
 (1) Recall from paragraph \ref{par:non-com} that the open subsets $\U_\F$ and $M_{\F,r}$ of~$M$ are defined by
\begin{eqnarray*}
    \U_\F&:=&\set{m\in M\mid \diff_m f_1\we\diff_m f_2\we\dots\we\diff_m f_s\neq0},\\
     M_{\F,r}&:=&\set{m\in M\mid \X_{f_1},\X_{f_2},\dots,\X_{f_r} \hbox{ are independent at } m}.
\end{eqnarray*}%
On $\U_\F\cap M_{\F,r}$ the functions $f_1,\dots,f_r$ define a foliation $\GG$ of corank $r$; similarly, the
functions $f_1,\dots,f_s$ define a foliation $\FF$ on it of rank $r$ (since $r+s=\dim M$). Obviously, $\FF$ is
contained in $\GG$. The condition that $\pb{f_i,f_j}=0$ for all $1\leqs i\leqs r$ and $1\leqs j\leqs s$, implies
that the Hamiltonian vector fields $\X_{f_1},\dots,\X_{f_r}$ are tangent to $\FF$ at every point. For all
$m\in\U_\F\cap M_{\F,r}$, the Hamiltonian vector fields $\X_{f_1},\dots,\X_{f_r}$ are independent at $m$, hence
they span $T_m\FF$. It follows that
\begin{eqnarray*}
  (T_m^\perp\GG)^{pol}&=&\set{\xi\in T_m^*M\mid \Pi_m(\xi,\diff_mf_i)=0\hbox{ for } i=1,\dots,r}\\
    &=&\set{\xi\in T_m^*M\mid \xi\left(\X_{f_i}(m)\right)=0\hbox{ for } i=1,\dots,r}=T_m^\perp\FF.
\end{eqnarray*}
It follows that $\FF$ is polar to $\GG$, hence $(\FF,\GG)$ is an abstract non-commutative integrable system.

\noindent
(2) Let $m$ be an arbitrary point of $M$. In a neighborhood $U$ of $m$, there exist smooth functions
    $f_1,\dots,f_s$, such that the level sets of $f_1,\dots,f_s$ and of $f_1,\dots,f_r$ define foliations, which
    coincide with $\FF$ and $\GG$ on $U$. Since $\FF$ is polar to $\GG$, the functions $f_1,\dots,f_r$ are in
    involution with the functions $f_1,\dots,f_s$ and the Hamiltonian vector fields of the functions $f_1,
    \dots,f_r$ are linearly independent at all points of $U$. It follows that $\F:=(f_1,\dots,f_s)$ is a
    non-commutative integrable system of rank $r$ on $U$.
\end{proof}

When both foliations $\FF$ and $\GG$ are given by fibrations $\F : M \to P $ and $\G : M \to L $ respectively, we
have a commutative diagram of submersive Poisson maps:
\begin{equation}\label{eq:Poissonsubm}
  \begin{diagram}
    \node{(M,\Pi)}\arrow{e,t}{\G}\arrow{s,l}{\F}
    \node{(L,\Pi_L)}\arrow{sw}\\
    \node{(P,\Pi_P)}
  \end{diagram}
\end{equation}%
where $\Pi_L$ is the zero Poisson structure on $L$. Moreover, item (3) in definition~ \ref{def:abstract-non-Com}
amounts for every $m \in M $ to the equality:
\begin{equation}\label{eq:duality}
  \F^* (T_{\F(m)}^*P) = \big(\G^*(T^*_{\G(m)} L)\big)^{pol}  .
\end{equation}
Conversely, it is clear that we have the following proposition:

\begin{prop}
  Suppose that (\ref{eq:Poissonsubm}) is a commutative diagram of submersive Poisson maps, where $L$ has dimension
  $r$ and is equipped with the zero Poisson structure $\Pi_L$, and $P $ has dimension $\dim M-r $.  If
  (\ref{eq:duality}) holds for every $m\in M$, then the pair of foliations $(\FF,\GG) $ defined on $M$ by $\F$ and
  $\G$ is an abstract non-commutative integrable system of rank $r$ on $(M,\Pi) $.
\end{prop}

% -*- mode: latex; tex-main-file: "art.tex"; -*-
%
\def\cprime{$'$}

\vfill

%\bibliographystyle{plain}
%\bibliography{ref}

\end{document}